\documentclass[11pt, reqno, a4paper]{amsart}
\usepackage{lmodern}
\usepackage[margin=1.2in]{geometry}
\numberwithin{equation}{section}
\usepackage{amssymb,amsfonts,amsthm}
\usepackage[utf8]{inputenc}
\usepackage{listings}
\usepackage{bm}
\usepackage[numbers,sort&compress]{natbib}
\usepackage[hyperpageref]{backref}
\usepackage{esint} 
\usepackage{color}
\usepackage{indentfirst}
\usepackage[bookmarks]{hyperref}
\usepackage{tikz}

\allowdisplaybreaks[4]

\newtheoremstyle{myremark}{10pt}{10pt}{}{}{\bfseries}{.}{.5em}{}

\newtheorem{remark}{Remark}[section]
 \newtheorem{thm}{Theorem}
 
 \newtheorem{lemma}{Lemma}[section]
 \newtheorem{prop}[thm]{Proposition}
 \theoremstyle{definition}

\usepackage{hyperref}

 \newcommand{\N}{\mathbb{N}}

\newcommand{\sn}[1]{[#1]_{W^{s,p,\alpha}(\mathbb{R}^{d})}}
\newcommand{\R}{\mathbb{R}}

\allowdisplaybreaks[4]

\begin{document}

\title[Failure of P\'olya--Szeg\H{o} inequality for weighted nonlocal seminorm]{Failure of P\'olya--Szeg\H{o} inequality for weighted nonlocal seminorm  and on nonlocal C-K-N inequality }

\author[P{.} Roy and V{.} Sahu]{Prosenjit Roy AND Vivek Sahu}

\address[P{.}R{.}]{Department of Mathematics and Statistics,
Indian Institute of Technology Kanpur, Kanpur - 208016, Uttar Pradesh, India}

\email{prosenjit@iitk.ac.in}

\address[V{.}S{.}]{Theoretical Statistics and Mathematics Unit,
Indian Statistical Institute, Delhi Centre, S.J. Sansanwal Marg, New Delhi, Delhi 110016, India}
\email{vivek@isid.ac.in,  viiveksahu@gmail.com}



\date{}

\dedicatory{}

\begin{abstract}
 We create counterexamples to show that P\'olya--Szeg\H{o} type of inequality fails for the weighted nonlocal seminorm 
\begin{equation*}
    [u]_{W^{s,p, \alpha}(\mathbb{R}^d)} := \left( \int_{\mathbb{R}^d} \int_{\mathbb{R}^d} \frac{|x|^{\alpha_{1} p} |y|^{\alpha_{2} p} |u(x)-u(y)|^{p}}{|x-y|^{d+sp}}  \, dx \, dy \right)^{\frac{1}{p}} \ , 
\end{equation*}
in all three regimes $\alpha_1 +\alpha_2 < 0, \, \alpha_1 +\alpha_2 = 0$, and $\alpha_1 +\alpha_2 >0$. These counterexamples are such that it completely settles the first two cases for all parameters. We also provide necessary and sufficient conditions on the parameters for $ [u]_{W^{s,p, \alpha}(\mathbb{R}^d)} < \infty$, thus  extending an earlier known result of \emph{[Serena Dipierro and Enrico Valdinoci, A density property for fractional weighted Sobolev spaces]}.
We also provide another proof of the fractional C-K-N type inequality, for most of the parameters (Hardy and Sobolev are included), using a familiar technique which is involved in the proof of fractional Sobolev inequality of  \emph{[Eleonora Di Nezza, Giampiero Palatucci, and Enrico Valdinoci, Hitchhiker's guide to the fractional
Sobolev spaces, Bulletin des sciences mathematiques]} together with a symmetrization technique. 
\end{abstract}

\maketitle


\section{Introduction}
Let $p\geq 1$, $s\in(0,1)$, and $\alpha,\alpha_1,\alpha_2\in\mathbb{R}$ satisfy $\alpha_1+\alpha_2=\alpha$. For a measurable set
$\Omega\subseteq\mathbb{R}^d$, we define the weighted fractional Gagliardo seminorm by
\begin{equation}\label{Defn: weighted Gagliardo seminorm}
    [u]_{W^{s,p, \alpha}(\Omega)} := \left( \int_{\Omega} \int_{\Omega} \frac{|x|^{\alpha_{1} p} |y|^{\alpha_{2} p} |u(x)-u(y)|^{p}}{|x-y|^{d+sp}} \, dx \, dy \right)^{\frac{1}{p}} \in [0, \infty], \quad \ \forall \ u \in L^{1}(\Omega) \ ,
\end{equation}
In the unweighted case, namely when $\Omega=\mathbb{R}^d$ and
$\alpha_1=\alpha_2=0$, the celebrated P\'olya--Szeg\H{o}-type inequality
for fractional Sobolev seminorms (see \cite[Section $9.2$]{Almgren1989} and \cite[Theorem A.$1.$, Appendix A]{frank2008non}) states that
$$ [u^*]_{W^{s,p}(\mathbb{R}^d)} \leq  [u]_{W^{s,p}(\mathbb{R}^d)},$$
where $u^*$ denotes the symmetric decreasing rearrangement
(Schwarz symmetrization) of $u$. In the above inequality, we have used the notation  $[u]_{W^{s,p, 0}(\mathbb{R}^d)} = [u]_{W^{s,p}(\mathbb{R}^d)}$, corresponding to the case  $\alpha_1 = \alpha_2=0$. In this article, we construct counterexamples to show that, in general, such an inequality does not hold for the weighted nonlocal seminorm $[u]_{W^{s,p, \alpha}(\mathbb{R}^d)}$ in any of the three cases $\alpha_1+\alpha_2<0$, $\alpha_1+\alpha_2=0$, and $\alpha_1+\alpha_2>0$.

\medskip

Before presenting the counterexamples, we claim that the range
$\alpha_1,\alpha_2\in(-d/p,s)$ is the non-trivial regime in which the
P\'olya--Szeg\H{o} inequality should be studied. In the complementary
region, the answer to the P\'olya--Szeg\H{o} inequality is trivial.
This is because, if $\alpha_1\leq -\frac{d}{p}$ or $\alpha_1\geq s$,
without any assumption on $\alpha_2$, then $[\varphi]_{W^{s,p,\alpha}(\mathbb{R}^d)}=\infty$, as shown in part (a) of Lemma \ref{lemma:critical interval} below.
Hence, the P\'olya--Szeg\H{o} inequality holds trivially in these cases. Dipierro and Valdinoci (see \cite[Lemma 2.1]{Serena2015}) showed that if $\alpha_1,\alpha_2\in(-d/p,s)$ and, in addition, $\alpha_1+\alpha_2>-d/p$, then $$[\varphi]_{W^{s,p,\alpha}(\mathbb{R}^d)}<\infty \qquad \text{for all } \varphi\in C_c^\infty(\mathbb{R}^d).$$ We believe that Lemma \ref{lemma:critical interval} is a new contribution to the literature, as it sharpens
\cite[Lemma 2.1]{Serena2015} by covering all the complementary cases.

\medskip

\begin{lemma}
    \label{lemma:critical interval}
    Let $u \in C_c^\infty(\mathbb{R}^d)$, $u  \not\equiv 0$. 
    \begin{itemize}
        \item[(a)] If $  \  \alpha_1 \leq -\frac{d}{p}$, then $[u]_{W^{s,p, \alpha}(\mathbb{R}^d)}= \infty$, or if $\alpha_1\geq s$, then $[u]_{W^{s,p, \alpha}(\mathbb{R}^d)}= \infty$.
        \item[(b)]  Assume $\alpha_1, \alpha_2 \in (-d/p,s)$, and define the \emph{order of vanishing of $u -u(0)$ at the origin} as
        \begin{equation*}
             m(u)= m:=\min\big\{|\beta|\ge1:\ \partial^\beta u(0)\ne0\big\}\in\N\cup\{\infty\}.
        \end{equation*}
 Then,   
 \begin{equation*}
     [u]_{W^{s,p, \alpha}(\mathbb{R}^d)}<\infty\quad\Longleftrightarrow\quad\alpha_1+\alpha_2>-\frac dp-(m-s).
 \end{equation*}
    \end{itemize}
\end{lemma}

\smallskip

We begin by presenting our first theorem, which provides a counterexample to the P\'olya--Szeg\H{o} inequality in the case $\alpha_1+\alpha_2<0$.

\begin{thm}\emph{[Counter example for $\alpha_1 + \alpha_2 <0$]}\label{thm:main}
Let $d\ge 1$, $p\ge 1$, $s\in(0,1)$, and let $\alpha_1,\alpha_2\in \mathbb{R}$ satisfy
\begin{equation*}
     -\frac{d}{p}<\alpha_1,\alpha_2<s
  \qquad\text{and}\qquad
  \alpha_1+\alpha_2<0 .
\end{equation*}
Then there is no constant $C>0$ such that
\begin{equation*}
[u^{*}]_{W^{s,p, \alpha}(\mathbb{R}^{d})} \leq C [u]_{W^{s,p, \alpha}(\mathbb{R}^{d})} \qquad\text{for all nonnegative } u \in C^{\infty}_{c}(\mathbb{R}^{d}). 
\end{equation*}
More precisely, let $\varphi\in C_c^\infty(\mathbb{R}^d)$ be radially nonincreasing with $0\le\varphi\le 1$, $\varphi\equiv 1$ on $B_1$ and $\operatorname{supp}\varphi\subset B_2$. Fix a unit vector $e$ and set $u_R(x)=\varphi(x-Re)$. Then,
\begin{enumerate}
  \item $u_R^*=\varphi$ for every $R>0$, and $0<[\varphi]_{W^{s,p, \alpha}(\mathbb{R}^{d})}<\infty$;
  \item there exists $C=C(d,p,s,\alpha_1,\alpha_2,\varphi)$ such that for all $R\ge 8$,
  \begin{equation*}
       [ u_{R}]_{W^{s,p, \alpha}(\mathbb{R}^{d})}^p\le C\left(R^{(\alpha_1+\alpha_2)p}+R^{(\alpha_1+\alpha_2-s)p}\right);
  \end{equation*}
  \item consequently $[ u_{R}]_{W^{s,p, \alpha}(\mathbb{R}^{d})}\to 0$ as $R\to\infty$.
\end{enumerate}
\end{thm}

Note that the above counterexample shows that the P\'olya--Szeg\H{o} inequality fails for all admissible choices of $d$, $p$, $s$, $\alpha_1$, and $\alpha_2$ satisfying the assumptions of the theorem. The reason for this failure is the interaction between translations and the negative total weight. Indeed, translating a function arbitrarily far from the origin does not change its Schwarz symmetrization. However, the weighted seminorm of the translated function tends to zero as the translation becomes large, while the weighted seminorm of its Schwarz symmetrization remains strictly positive. 

\smallskip

We next consider the case $\alpha_1+\alpha_2=0$ and show that the P\'olya--Szeg\H{o} inequality also fails in this case. Unlike the previous theorem, the counterexample here is based on a strict difference between the weighted seminorms of a function and its Schwarz symmetrization.

\begin{thm}\emph{[Counter example for $\alpha_1 + \alpha_2 =0$]}\label{thm:zero}
Let $d\ge1$, $p\ge1$, $s\in(0,1)$, and let $\alpha_1=\beta$, $\alpha_2=-\beta$ with $0<|\beta|<\min\Big\{s,\frac dp\Big\}.$
Let $\varphi$ and $u_R$ be as in Theorem~\ref{thm:main}, and $\delta_{\beta}$ be as defined in \eqref{Defn: delta_beta}. Then,
\begin{enumerate}
  \item $\delta_\beta(\varphi)>0$; in particular $\sn{\varphi}>[\varphi]_{W^{s,p}(\R^d)}$;
  \item $\sn{u_R}^p\to[\varphi]^p_{W^{s,p}(\R^d)}$ as $R\to\infty$;
  \item consequently
  \[
    \lim_{R\to\infty}\Big(\sn{u_R^*}^p-\sn{u_R}^p\Big)=\delta_\beta(\varphi)>0,
  \]
  and there exists $R_0$ such that $\sn{u_R}<\sn{u_R^*}$ for all $R\ge R_0$.
\end{enumerate}
In particular, the P\'olya--Szeg\H{o} inequality $\sn{u^*}\le\sn{u}$ fails for nonnegative $u\in C_c^\infty(\R^d)$.
\end{thm}

\smallskip

Notice that for $d\ge 2$ and $\alpha_1,\alpha_2\ge 0$ (with $\varphi$ radial) the translation of
a radial function cannot decrease the seminorm. Define
\[
f(h):=[\varphi(\cdot-h)]^p_{W^{s,p,\alpha}(\mathbb{R}^{d})}=\int_{\mathbb{R}^{d}} \int_{\mathbb{R}^{d}} G_{x,y}(h)\,d\mu_\varphi(x,y),\qquad
G_{x,y}(h):=|x+h|^{\alpha_{1}p}|y+h|^{\alpha_{2}p},
\]
where $d\mu_\varphi(x,y)=\frac{|\varphi(x)-\varphi(y)|^p}{|x-y|^{d+sp}}\,dx\,dy$.
Since $\varphi$ and the weights are radial, $f$ is radial, so averaging over $|h'|=|h|$ gives
\[
f(h)=\int_{\mathbb{R}^{d}} \int_{\mathbb{R}^{d}} \Big(\fint_{|h'|=|h|}G_{x,y}(h')\,d\sigma(h')\Big)\,d\mu_\varphi(x,y).
\]
For $d\ge 2$ the function $h\mapsto\log G_{x,y}(h)=\alpha_{1}p \log|x+h|+\alpha_{2}p\log|y+h|$ is subharmonic, because
$\log|\cdot|$ is subharmonic and $\alpha_{1}p,\alpha_{2}p\ge0$; hence $G_{x,y}=\exp(\log G_{x,y})$ is subharmonic, being
the exponential (convex, nondecreasing) of a subharmonic function. So its spherical mean dominates $G_{x,y}(0)$, and
therefore $f(h)\ge f(0)$. (A product of nonnegative subharmonic functions need not be subharmonic;
the logarithm is what makes the argument work.) Thus translations of radial functions cannot
produce counterexamples when $d\ge2$ and $\alpha_1,\alpha_2\ge0$.
For $d=1$ this argument breaks down: when $\alpha_1=0$, $\alpha_2p=b\in(0,1)$, one has
$G_{x,y}(h)=|y+h|^{b}$ and $G_{x,y}''(h)=b(b-1)|y+h|^{b-2}<0$ for $h\neq -y$, so $G_{x,y}$ is
concave away from its cusp at $h=-y$. This is exactly the leverage used in Theorem \ref{thm:onesided}.

\smallskip 

For $\varphi\in C_c^\infty(\R^d)$, define its \emph{energy marginal} by
\[
  m_\varphi(\eta):=\int_{\R^d}\frac{|\varphi(\xi)-\varphi(\eta)|^p}{|\xi-\eta|^{d+sp}}\,d\xi .
\] Let $\varphi(\cdot-h)$ denote the translation of $\varphi$ by $h \in \mathbb{R}^d.$ For simplicity, in the case $\alpha_{1}=0$ and $\alpha_{2}= \varepsilon>0$, we use the following notation for weighted Gagliardo seminorm,
\begin{equation}\label{Defn: W0,epsilon}
    [u]^{p}_{0,\varepsilon} := \int_{\mathbb{R}^{d}} \int_{\mathbb{R}^{d}} \frac{|y|^{\varepsilon p}|u(x)-u(y)|^{p}}{|x-y|^{d+sp}} dx \, dy.
\end{equation}

\smallskip

\begin{thm}\emph{[Counter example for $\alpha_1 + \alpha_2 > 0$]}\label{thm:onesided}
Let $d=1$,  $\alpha_1=0$ and $\alpha_2=\varepsilon>0$, $p\ge1$, $s\in(0,1)$ with $sp\ge1$, and set $\varepsilon_0:=\min\{s/2,\,1/(2p)\}$. There is an even function $\varphi\in C_c^\infty(\R)$, nonincreasing on $[0,\infty)$, such that for every $\varepsilon\in(0,\varepsilon_0]$ there exists $h=h(\varepsilon)>0$ with
\[
  [\varphi(\cdot-h)]_{0,\varepsilon}\ <\ [\varphi]_{0,\varepsilon}=[\varphi(\cdot-h)^*]_{0,\varepsilon}.
\]
More precisely, for any such $\varphi$ with $\varphi\equiv1$ on $[-1,1]$ and any $b=\varepsilon p\in(0,1)$, $\varepsilon<s$,
\[
  \lim_{h\to0^+}\frac{[\varphi(\cdot-h)]_{0,\varepsilon}^p-[\varphi]_{0,\varepsilon}^p}{h^2}=\frac{b(b-1)}{2}\,Q_b(\varphi),
  \quad \text{where }
  Q_b(\varphi):=\int_{\R}|\eta|^{b-2}\big(m_\varphi(\eta)-m_\varphi(0)\big)\,d\eta ,
\]
and the integral defining $Q_b(\varphi)$ converges absolutely. Hence the inequality fails whenever $Q_b(\varphi)>0$.
\end{thm}

\medskip
\noindent\textbf{Scope of Theorem \ref{thm:onesided} and open cases:}
Theorem \ref{thm:onesided} is a counterexample only in the following setting: $d=1$, $\alpha_1=0<\alpha_2=\varepsilon$ (the case $\alpha_2=0<\alpha_1$ follows by symmetry of the seminorm), $sp\ge 1$ and $0<\varepsilon\le\varepsilon_0$. For $\alpha_1+\alpha_2>0$ the validity of the
P\'olya--Szeg\H{o} inequality is therefore open in all other situations, in particular:
\begin{itemize}
\item $d=1$, $\alpha_1=0$ and $sp<1$;
\item $d=1$, $\alpha_1=0$, $sp\ge 1$ and $\varepsilon>\varepsilon_0$;
\item $d=1$ with $\alpha_1\neq 0$ and $\alpha_2\neq 0$;
\item $\alpha_1<0<\alpha_2$ with $\alpha_1+\alpha_2>0$, in any dimension;
\item $d\ge 2$. Here we only know that translations of radial functions cannot
      produce counterexamples when $\alpha_1,\alpha_2\ge 0$ (see the discussion after
      Theorem \ref{thm:zero}), and that the inequality holds for $p=1$, $\alpha_1=0$ for functions whose
      superlevel sets are balls (Proposition \ref{prop:partial-positive}).
\end{itemize}
Moreover, Theorems \ref{thm:zero} and \ref{thm:onesided}  show the failure of the inequality with constant $1$, whereas Theorem \ref{thm:main}
 shows that no constant $C$ works.


\smallskip

The natural question that arises is to determine the regime in which the P\'olya--Szeg\H{o} inequality holds. We establish partial results for $\alpha_1=0$, $\alpha_2>0$, and $p=1$, under the assumption that the level sets of $u$ are balls.

\smallskip 

\noindent\textbf{The weighted fractional perimeter:}
Throughout the rest of this section, we take $p=1$ and write
$\alpha_2=b=\varepsilon$. For a measurable set $E\subset\mathbb{R}^d$
with $|E|<\infty$, define
\begin{equation}\label{Defn: Pb}
    P_b(E):=[\mathbf 1_E]_{0,\varepsilon}:=\int_{\mathbb{R}^d}\int_{\mathbb{R}^d}
  \frac{\bigl(|x|^0|y|^b\bigr)
  |\mathbf{1}_E(x)-\mathbf{1}_E(y)|}
  {|x-y|^{d+s}}\,dx\,dy\in[0,\infty].
\end{equation}
A simple calculation yields
\[
  P_{b}(E)=\int_E\int_{E^c}\frac{|x|^b+|y|^b}{|x-y|^{d+s}}\,dy\,dx ,
\]
which represents the interaction between $E$ and its complement, with each pair $(x,y)$ weighted by $|x|^b+|y|^b$. When $b=0$, this reduces to twice the usual fractional $s$-perimeter, 
$$P^s(E)=\int_E\int_{E^c}|x-y|^{-d-s}\,dy\,dx.$$

By the coarea formula, proving the P\'olya--Szeg\H{o} inequality in this case (see Lemma \ref{cor:sets}) is equivalent to proving the isoperimetric inequality $P_b(E) \geq P_b(E^*)$, where $E^*$ is the ball centred at the origin satisfying $|E|=|E^*|$.

\smallskip 

\begin{prop}\emph{\textbf{[A partial positive result]}}
\label{prop:partial-positive}
Let $d\ge2$, $p=1$ and $\alpha_1 = 0, 0 < \alpha_2=:\varepsilon $. Let $u\ge0$ be measurable and vanishes at infinity, and suppose that for a.e.\ $t>0$ the superlevel set $\{u>t\}$ is (up to a null set) a ball $B_{\rho(t)}(h(t))$, where the centres $h(t)\in\R^d$ may depend on $t$. Then
\[
  [u^*]_{0,\varepsilon}\le[u]_{0,\varepsilon}.
\]
\end{prop}


The main idea of the proof is the same as that in the discussion following the statement of Theorem \ref{thm:zero}, where translations of radial functions are ruled out as counterexamples when $d\geq2$. Note that translations of radial functions indeed satisfy the assumptions of our result.

\subsection{On nonlocal Caffarelli, Kohn and Nirenberg Inequality} Caffarelli, Kohn and Nirenberg in  \cite{Caffarelli1984} established the so-called ``Caffarelli-Kohn-Nirenberg inequality" or ``C-K-N inequality" in short. In  \cite[page no. 1]{Caffarelli1984}, they proved the following theorem: let  $p,  ~ q \geq 1,  ~\tau >0,  ~ a \in (0,1], ~  \alpha,  ~ \beta, ~ \gamma \in \mathbb{R}$ such that
\begin{equation*}
    \frac{1}{\tau} + \frac{\gamma}{d}, \quad \frac{1}{p}+ \frac{\alpha}{d}, \quad \frac{1}{q} + \frac{\beta}{d} >0,
\end{equation*}
\begin{equation*}
    \frac{1}{\tau} + \frac{\gamma}{d} = a \left( \frac{1}{p} + \frac{\alpha-1}{d} \right) + (1-a)\left( \frac{1}{q} + \frac{\beta}{d} \right),
\end{equation*}
and, with  $\gamma = a \sigma + (1-a)\beta$,
\begin{eqnarray*}
    0 \leq \alpha - \sigma
\end{eqnarray*}
and
\begin{equation*}
    \alpha -\sigma \leq 1 \quad \text{if} \quad \frac{1}{\tau} + \frac{\gamma}{d} = \frac{1}{p} + \frac{\alpha-1}{d}.
\end{equation*}
Then, for every  $u \in C^{1}_{c}(\mathbb{R}^d)$,
\begin{equation}\label{ckn}
    \| |x|^{\gamma} u \|_{L^{\tau}(\mathbb{R}^d)} \leq C \| |x|^{\alpha} \nabla u \|^{a}_{L^{p}(\mathbb{R}^d)}\| |x|^{\beta} u \|^{(1-a)}_{L^{q}(\mathbb{R}^d)},
\end{equation}
for some positive constant  $C$ independent of  $u$. The C-K-N inequality has been a topic of research interest since its inception. We refer  to some of the work of  \cite{abdellaoui2005some, bartsch2007existence, catrina2000caffarelli,  catrina2001caffarelli,  chern2010minimizers, dolbeaulty2009symmetry,   kristaly2013caffarelli, radulescu2011, nguyen2019hardy, wang2003caffarelli,   wei2022stability} and references therein. 

\smallskip

Nguyen and Squassina in  \cite{Nguyen2018}, extended the local C-K-N  inequalities to fractional Sobolev spaces. They considered the following weighted Gagliardo seminorm $[u]_{W^{s,p}(\Omega)}$
and proved the following theorem: Let  $d\geq 1$,  $p>1$,  $q \geq 1$,  $\tau >0$,  $a \in (0,1]$,  $\alpha, ~\beta, 
~\gamma \in \mathbb{R}$ be such that
\begin{equation*}
    \frac{1}{\tau}+ \frac{\gamma}{d} = a \left( \frac{1}{p} + \frac{\alpha-s}{d} \right) + (1-a) \left( \frac{1}{q} + \frac{\beta}{d} \right)  \ .
\end{equation*}
In the case  $a>0$, assume in addition that, with  $\gamma = a \sigma + (1-a) \beta$,
\begin{equation*}
    0 \leq \alpha - \sigma
\end{equation*}
and 
\begin{equation*}
    \alpha - \sigma \leq s \hspace{.5cm} \text{if} \hspace{.5cm} \frac{1}{\tau}+ \frac{\gamma}{d} = \frac{1}{p}+ \frac{\alpha-s}{d}  .
\end{equation*}
If  $\frac{1}{\tau} + \frac{\gamma}{d} >0$, then
      \begin{equation}\label{nguyen}
          \| |x|^{\gamma}u \|_{L^{\tau}(\mathbb{R}^d)} \leq C [u]^{a}_{W^{s,p, \alpha}(\mathbb{R}^d)} \||x|^{\beta} u \|^{(1-a)}_{L^{q} (\mathbb{R}^d)},  \quad  \forall \ u  \in C^{1}_{c}(\mathbb{R}^d)  .
      \end{equation}
      
\smallskip

Clearly,  \eqref{nguyen} is the appropriate fractional analog of  \eqref{ckn}.  The method of proof used by Nguyen and Squassina in  \cite[Theorem 1.1 and Theorem 3.1]{Nguyen2018} involves a rigorous use of Sobolev inequality for  $W^{s,p}$ functions, a clever use of summation technique, and appropriate sharp estimates. The main aim of this article is to present an alternative proof of  \eqref{nguyen} (a fractional Caffarelli-Kohn-Nirenberg type inequality)  under some additional assumptions. More precisely, under the assumptions  $ \gamma \leq 0, ~\alpha_1 = \alpha_2=0$ and  $\beta =0$.  In the class of nonnegative, radially nonincreasing functions, these additional restrictions ($\alpha_1 = \alpha_2=0$ and  $\beta =0$) can be removed.  Our approach is inspired by the technique in the proof of fractional Sobolev inequality  \cite[Theorem 6.5]{di2012hitchhikers} in the pioneering work of  Di Nezza, Palatucci, and Valdinoci    \cite{di2012hitchhikers} (commonly known as Hitchhiker's guide).  This together with a simple symmetrization argument is sufficient to establish our main theorem.  In view of Theorems \ref{thm:main}, \ref{thm:zero}, and \ref{thm:onesided}, the P\'olya--Szeg\H{o} inequality does not hold in general for weighted nonlocal seminorms. Therefore, we impose the additional assumption that $u$ is nonnegative, radially nonincreasing function. More precisely, we state our main result below.


\begin{thm}\label{Main result}
Let    $p , ~q \geq 1, ~ \tau \geq \max\{p,q\}, ~ \alpha \in \mathbb{R}$,   $s \in (0,1)$,  $\gamma \leq 0$,   $a \in (0,1]$, and  $ \beta \in \mathbb{R}$ be such that $\alpha=\alpha_1+\alpha_2$, $-\frac{d}{p}<\alpha_1 , \alpha_{2}<s$, and
     \begin{equation*}
 \frac{1}{\tau} + \frac{\gamma}{d}, \hspace{3mm}    \frac{1}{p} +\frac{\alpha - s}{d}, \hspace{3mm} \frac{1}{q} + \frac{\beta}{d} >0,
\end{equation*}
and
\begin{equation*}
        \frac{1}{\tau} + \frac{\gamma}{d} =  a \left( \frac{1}{p} + \frac{\alpha-s}{d} \right) + (1-a) \left( \frac{1}{q} + \frac{\beta}{d} \right).
    \end{equation*} Then there exists a constant  $C>0$ (independent of  $u$) such that
    \begin{equation*}
        \| |x|^{\gamma}u \|_{L^{\tau}(\mathbb{R}^d)} \leq C [u]^{a}_{W^{s,p, \alpha}(\mathbb{R}^d)} \| |x|^{\beta} u \|^{(1-a)}_{L^{q}(\mathbb{R}^d)}, \  \hspace{.3cm}   \text{for all radially nonincreasing } \, u  \in C^{1}_{c}(\mathbb{R}^d).
    \end{equation*}
     Moreover, when  $\alpha_{1}=\alpha_{2}=0$ and  $\beta = 0$, the assumption of  $u$ to be radially nonincreasing function in the previous inequality can be removed.
\end{thm}

\smallskip

Let  $\gamma=-s, ~ \tau =p, ~ \alpha_1=\alpha_2=0$, and $a=1$ in the above theorem. Then we obtain
\begin{equation}\label{fractionalhardy}
    \int_{\mathbb{R}^{d}} \frac{|u(x)|^{p}}{|x|^{sp}} dx \leq C \int_{\mathbb{R}^d}  \int_{\mathbb{R}^d} \frac{|u(x)-u(y)|^{p}}{|x-y|^{d+sp}} \, dx \, dy =: [u]^p_{W^{s,p}(\mathbb{R}^d)} , \hspace{.3cm} \forall \ u \in C^{1}_{c}(\mathbb{R}^d),
\end{equation}
with $sp <d$, which is the well-known ``fractional Hardy inequality". For some well-known proofs of \eqref{fractionalhardy}, we refer to  \cite{frank2008non, Mazya2002}. The proof of the above inequality also follows directly from a simple modification of the proof of the fractional Sobolev inequality (see \cite[Theorem $6.5$]{di2012hitchhikers} for the proof of the fractional Sobolev inequality). The main steps are roughly as follows:

\medskip

\noindent \textbf{Short proof of fractional Hardy inequality: } In  \cite{di2012hitchhikers}, it is proved (as a combination of \cite[Lemma $6.2$]{di2012hitchhikers} and \cite[Lemma $6.3$]{di2012hitchhikers}) that for  $sp<d$,   $k \in \mathbb{Z}$,  $A_k := \{x: |u(x)| > 2^k\}$,  $a_k:= |A_k|$,  $D_k:= A_k\setminus A_{k+1}$, and  $d_{k}:= |D_{k}|$, one has, for some constant  $C>0$ (independent of $u$), 
\begin{equation}\label{hardy1}
    [u]^{p}_{W^{s,p}(\mathbb{R}^d)} \geq C\sum_{k \in \mathbb{Z}} 2^{k p} d^{ \frac{d-sp}{d}}_{k}.
\end{equation}
Now, using upper estimates of  $u$ on  $D_k$ together with a symmetrization argument and the condition  $sp<d$, we have
\begin{align}\label{hardy2}
    \int_{\mathbb{R}^d} \frac{|u(x)|^{p}}{|x|^{sp}} dx \leq 2^{p} \sum_{k \in \mathbb{Z}} 2^{kp} \int_{D_{k}} \frac{1}{|x|^{sp}} dx & \leq 2^{p} \sum_{k \in \mathbb{Z}} 2^{kp} \int_{D_{k}^*} \frac{1}{|x|^{sp}} dx \nonumber \\ & = \frac{2^{p} d^{\frac{d-sp}{d}}}{(d-sp)| \mathbb{S}^{d-1}|^{\frac{-sp}{d}}}  \sum_{k \in \mathbb{Z}} 2^{kp} d^{\frac{d-sp}{d}}_{k},
\end{align}
where  $D_k^*$ denotes the ball centered at origin such that  $|D_k|=|D_k^*|$, and $\mathbb{S}^{d-1}$ represents the boundary of an open ball of radius  $1$ centered at  $0$ in  $\mathbb{R}^{d}$. Therefore, combining  \eqref{hardy1} and  \eqref{hardy2} gives  \eqref{fractionalhardy} for  $sp<d$.

\smallskip


\smallskip 

The available literature is huge on Hardy inequalities and fractional Hardy inequalities, and we refer to some of the related work in  \cite{adi2002,  Brezis1997, Brezis2000}, and for the generalization of Hardy-type inequalities to fractional Sobolev spaces, we refer to some of the work in  \cite{ADICCM2026, adiarka2018, AdiJFA2026, AdiCVPDE2026,  Adi2025, BanerjeeGangulySahu, bogdan2011dyda, Gyula2026, dyda2004fractional, dyda2015, dyda2024, Dyda2024JFA, frank2008non, KijaczkoSNS, VivekCCM2026, VivekQS126, loss2010hardy, mallick2022caffarelli, VivekJMAA2025, VivekQS26} references therein. 

For more literature on fractional C-K-N inequalities, we refer to the work of \cite{Abdellaoui2017, Ao2022, Chen2021, DeNitti2025, Gonzalez2025, rawat2026}.

\smallskip

\textbf{The article is organised as follows.} In Section \ref{Sec2}, we establish necessary and sufficient conditions for the finiteness of the weighted fractional Gagliardo seminorm of smooth, compactly supported functions, and prove Lemma \ref{lemma:critical interval}. In Section \ref{sec3}, we investigate the validity of the P\'olya--Szeg\H{o} inequality for weighted fractional Gagliardo seminorms. In particular, we prove Theorems \ref{thm:main}, \ref{thm:zero} and \ref{thm:onesided}, establishing counterexamples in the cases $\alpha_1+\alpha_2<0$ and $\alpha_1+\alpha_2=0$, as well as in a particular case with $\alpha_1+\alpha_2>0$. We also establish a partial positive result for $p=1$ under an additional assumption on the superlevel sets of the function. Finally, in Section \ref{proof of main result}, we prove the main result, Theorem~\ref{Main result}, which establishes a fractional Caffarelli--Kohn--Nirenberg-type inequality.

\section{Proof of Lemma \ref{lemma:critical interval}}\label{Sec2}

In this section, we prove Lemma \ref{lemma:critical interval}, which
provides necessary and sufficient conditions for the finiteness of the
weighted fractional Gagliardo seminorm of a smooth, compactly supported
function. We first identify the parameter ranges in which the seminorm
is necessarily infinite. We then establish a sharp criterion for its
finiteness in terms of the sum $\alpha_1+\alpha_2$ and the order of
vanishing of the function at the origin.

\begin{proof}[\textbf{Proof of Lemma \ref{lemma:critical interval}}] \textbf{Proof of (a) :\, } Let us denote $K(z):= |z|^{-(d+sp)}, \ a= \alpha_1p$ and $c=\alpha_2p$. Since $\{|u|>0\}=\bigcup_{n}\{|u|>1/n\}$ has positive measure, there is $t>0$ such that
\[
  A:=\{|u|>2t\}\quad\text{satisfies}\quad 0<|A|<\infty .
\]
As $|A\cap B_R|\to|A|$ for $R\to\infty$ and $|A\cap B_\delta|\to0$ for $\delta\to0$, there are $0<\delta<R$ with
\[
  A':=A\cap(B_R\setminus B_\delta),\qquad |A'|>0,\qquad \delta\le|y|\le R\ \text{ on }A'.
\]
Let $C:=\{|u|<t/2\}$. Its complement $\{|u|\ge t/2\}$ has finite measure since $u$ vanishes at infinity. The annulus $D_\rho:=B_{2\rho}\setminus B_\rho$ has $|D_\rho|=(2^d-1)|B_1|\rho^d\to\infty$, so there is $\rho\ge\max\{R,1\}$ with
\[
  C':=C\cap D_\rho,\qquad |C'|\ge|D_\rho|-|\R^d\setminus C|>0,\qquad \rho\le|y|\le2\rho\ \text{ on }C'.
\]
We shall use the elementary pointwise bounds
\begin{equation}\label{eq:gaps}
  |u(x)-u(y)|>
  \begin{cases}
    t & \text{if } |u(x)|\le t,\ y\in A,\\
    t/2 & \text{if } |u(x)|>t,\ y\in C,\\
    t & \text{if } x\in A,\ y\in C,
  \end{cases}
\end{equation}
which follow from $|u(x)-u(y)|\ge\big||u(x)|-|u(y)|\big|$.

\medskip
\emph{Proof of  the first part $\alpha_1 \leq - d/p$ of (a):}  \ Assume $a=\alpha_1p\le-d$. By Tonelli's theorem,
\[
  [u]^p_{W^{s,p,\alpha}(\mathbb{R}^{d})}=\int_{\R^d}|x|^{a}\,I(x)\,dx,\qquad I(x):=\int_{\R^d}|y|^{c}\,|u(x)-u(y)|^p\,K(x-y)\,dy .
\]
We claim that $I(x)\ge\kappa>0$ for a.e.\ $x\in B_1$. Fix such an $x$.

\smallskip
\noindent\emph{(i) If $|u(x)|\le t$}, we integrate only over $y\in A'$. There $|u(x)-u(y)|>t$ by \eqref{eq:gaps}, $|x-y|\le1+R$ so $K(x-y)\ge(1+R)^{-d-sp}$, and $|y|^c\ge m_A:=\min\{\delta^c,R^c\}>0$. Hence
\[
  I(x)\ \ge\ t^p\,m_A\,(1+R)^{-d-sp}\,|A'|=:\kappa_1>0 .
\]

\noindent\emph{(ii) If $|u(x)|>t$}, we integrate only over $y\in C'$. There $|u(x)-u(y)|>t/2$ by \eqref{eq:gaps}, $|x-y|\le1+2\rho$, and $|y|^c\ge m_C:=\min\{\rho^c,(2\rho)^c\}>0$. Hence
\[
  I(x)\ \ge\ (t/2)^p\,m_C\,(1+2\rho)^{-d-sp}\,|C'|=:\kappa_2>0 .
\]
Thus $I(x)\ge\kappa:=\min\{\kappa_1,\kappa_2\}$ for a.e.\ $x\in B_1$, and since $a+d\le0$,
\[
  [u]^p_{W^{s,p,\alpha}(\mathbb{R}^{d})}\ \ge\ \kappa\int_{B_1}|x|^{a}\,dx=\kappa\,|S^{d-1}|\int_0^1r^{a+d-1}\,dr=\infty .
\]

\medskip
\emph{Proof of the second part $\alpha_2  \geq s$ of (a):} \  Notice by symmetry doing for $\alpha_2$  (instead of $\alpha_1$) also suffices.  Assume $\alpha_2\ge s$, i.e.\ $c\ge sp$. Restricting the integration to $x\in A'$ and $y\in C\cap\{|y|\ge2R\}$,
\[
  [u]^p_{W^{s,p,\alpha}(\mathbb{R}^{d})}\ \ge\ \int_{A'}|x|^{a}\Big(\int_{C\cap\{|y|\ge2R\}}|y|^{c}\,|u(x)-u(y)|^p\,K(x-y)\,dy\Big)dx .
\]
For such $x,y$ we have $|u(x)-u(y)|>t$ by \eqref{eq:gaps}, and $|x-y|\le R+|y|\le2|y|$, so $K(x-y)\ge2^{-d-sp}|y|^{-d-sp}$. Hence the inner integral is greater than or equal to
\[
  t^p\,2^{-d-sp}\int_{C\cap\{|y|\ge2R\}}|y|^{\,c-d-sp}\,dy.
\]
 Since $u \in C_c^\infty(\mathbb{R}^d),$ there exist a large $M$ such that $B^{c}_M \subset C\cap\{|y|\ge2R\} $. Then  clearly $\int_{\{|y|\ge M\}}|y|^{c-d-sp}\,dy=\infty$ because $c-d-sp\ge-d$.  This proves (a).
   
\bigskip

\textbf{Proof of (b) :\ } Set $\gamma:=a+c+d+p\,(m-s)(\gamma=+\infty\text{ if }m=\infty)$ and $ F_u(x,y):=|x|^{a}\,|y|^{c}\,\frac{|u(x)-u(y)|^p}{|x-y|^{d+sp}}$.
Fix $L\ge1$ with $\operatorname{supp}u\subset B_L$,  and put $M:=\|u\|_\infty$, $\Lambda:=\|\nabla u\|_\infty$. We split
\[
  \R^d\times\R^d=Q\ \cup\ Q^c,\qquad Q:=B_1\times B_1 .
\]
Step 1 shows that $\int_{Q^c} \int_{Q^c} F_u<\infty$ whenever $a,c\in(-d,sp)$, regardless of $a+c$ and of $u$ near $0$. Steps~2--4 show that $\int_Q \int_Q F_u<\infty$ if and only if $\gamma>0$.

\medskip
\emph{Step 1: On the region $Q^c$ .} The set $Q^c$ consists of three kinds of pairs: $\{|x|\ge1,\ |y|\ge1\}$ and the two \emph{mixed} regions $\{|x|\ge1,\ |y|<1\}$ and $\{|x|<1,\ |y|\ge1\}$. All three are covered by
\[
  Q^c=\{|x|\ge1\}\times\R^d\ \cup\ \R^d\times\{|y|\ge1\},\qquad\text{so}\qquad \mathbf 1_{Q^c}\le\mathbf 1_{\{|x|\ge1\}}+\mathbf 1_{\{|y|\ge1\}},
\]
where in the first set $y$ ranges over \emph{all} of $\R^d$ (in particular over $B_1$, which gives the mixed region $\{|x|\ge1,|y|<1\}$), and symmetrically in the second. Hence
\[
  \int_{Q^c} \int_{Q^c} F_u\ \le\ \int_{|x|\ge1}|x|^{a}I(x)\,dx+\int_{|y|\ge1}|y|^{c}J(y)\,dy,
\]
with $I$ as below and $J(y):=\int|x|^a|u(x)-u(y)|^p|x-y|^{-d-sp}\,dx$. By symmetry (exchange $(x,a)$ and $(y,c)$) it suffices to show that
\[
  \int_{|x|\ge1}|x|^{a}\,I(x)\,dx<\infty,\qquad I(x):=\int_{\R^d}|y|^{c}\,\frac{|u(x)-u(y)|^p}{|x-y|^{d+sp}}\,dy .
\]
We use $|u(x)-u(y)|\le\min\{\Lambda|x-y|,\,2M\}$, and that $u(x)-u(y)=0$ if $|x|,|y|\ge L$. Fix $x$ with $|x|\ge1$ and split the $y$-integral.

(i) \emph{$|y-x|<|x|/2$.} Then $\frac12|x|<|y|<\frac32|x|$, so $|y|^c\le C|x|^c$; the integrand vanishes unless $|x|\le2L$ (otherwise $|x|,|y|>L$). Hence this part is at most
\[
  C\,\Lambda^p\,|x|^{c}\,\mathbf 1_{\{|x|\le2L\}}\int_{|z|<|x|/2}|z|^{p-d-sp}\,dz\ \le\ C\,\mathbf 1_{\{|x|\le2L\}},
\]
because $p-d-sp>-d$.

(ii) \emph{$|y-x|\ge|x|/2$ and $|x|\le L$.} This is where the mixed region $\{|x|\ge1,\ |y|<1\}$ sits: $y$ may approach the origin, where $|y|^c$ is singular, but it stays at distance $\ge\frac12$ from $x$. Thus $|x-y|\ge\frac12$ and this part is at most
\[
  (2M)^p\int_{\R^d}|y|^{c}\,\min\{2^{d+sp},\,|x-y|^{-d-sp}\}\,dy\ \le\ C,
\]
uniformly in $|x|\le L$, because $|y|^c$ is locally integrable ($c>-d$) and $|y|^{c-d-sp}$ is integrable at infinity ($c<sp$).

(iii) \emph{$|y-x|\ge|x|/2$ and $|x|>L$.} The integrand vanishes unless $|y|\le L$ (this again includes $y$ near the origin), and then $|x-y|\ge|x|/2$, so this part is at most
\[
  (2M)^p\,2^{d+sp}|x|^{-d-sp}\int_{B_L}|y|^{c}\,dy=C\,|x|^{-d-sp}.
\]
Therefore $I(x)\le C\big(\mathbf 1_{\{|x|\le2L\}}+|x|^{-d-sp}\big)$ for $|x|\ge1$, and
\[
  \int_{|x|\ge1}|x|^{a}I(x)\,dx\ \le\ C\int_{1\le|x|\le2L}|x|^a\,dx+C\int_{|x|\ge1}|x|^{a-d-sp}\,dx<\infty,
\]
since $a<sp$. The same argument with the roles of $x$ and $y$ exchanged handles $\{|y|\ge1\}$.

\medskip
\emph{Step 2: dyadic decomposition of $Q$ and scaling.} For $(x,y)\in Q\setminus\{(0,0)\}$ let $\rho(x,y):=\max\{|x|,|y|\}$ and
\[
  S_k:=\{(x,y):\ 2^{-k-1}<\rho(x,y)\le2^{-k}\},\qquad k=0,1,2,\dots,
\]
so that $Q\setminus\{(0,0)\}$ is the disjoint union of the shells $S_k$, and $S_k=2^{-k}S_0$ (see Figure~\ref{fig:shells}).
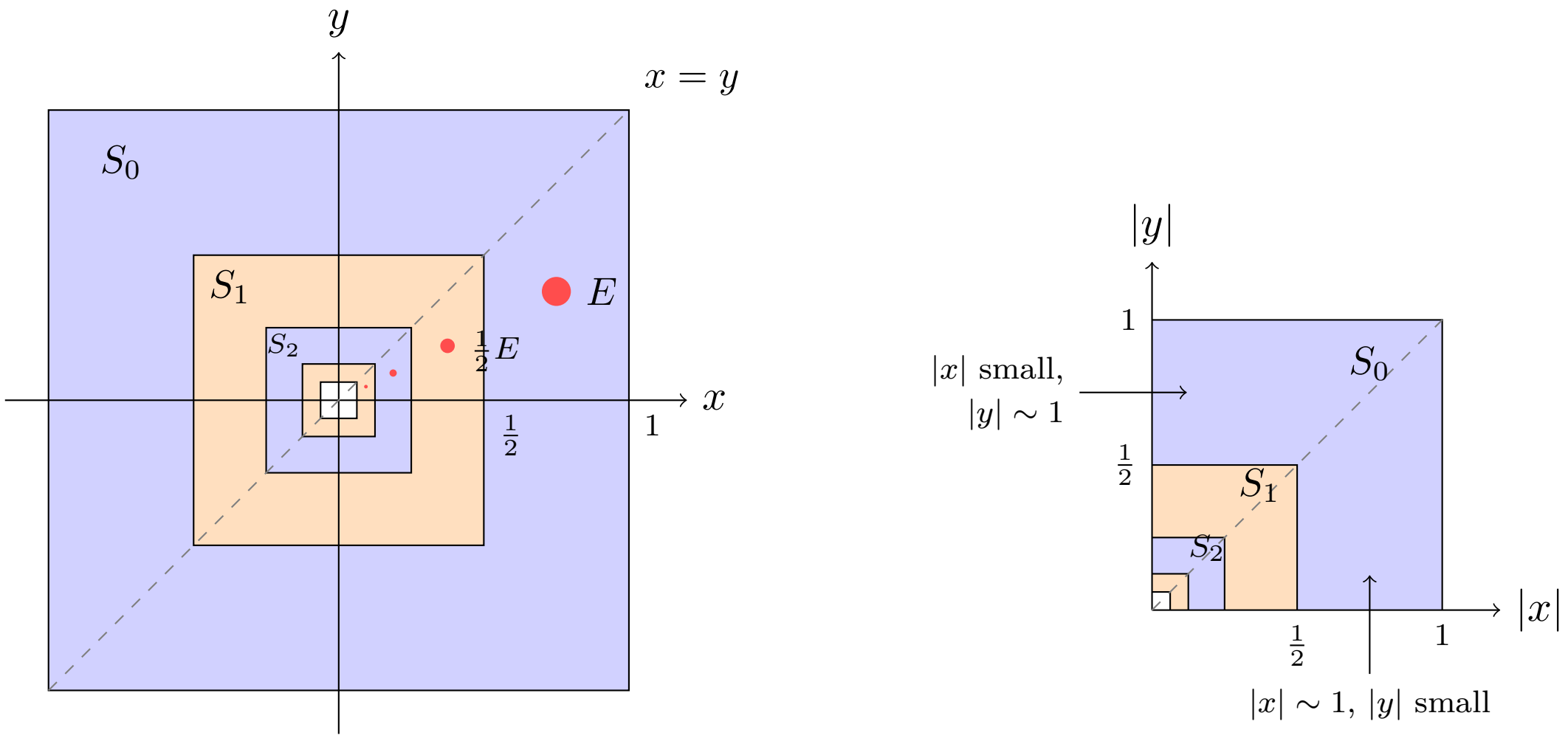
\begin{figure}[ht]
\centering
\begin{tikzpicture}[scale=2.6]
  \fill[blue!18]  (-1,-1) rectangle (1,1);
  \fill[white]    (-0.5,-0.5) rectangle (0.5,0.5);
  \fill[orange!25](-0.5,-0.5) rectangle (0.5,0.5);
  \fill[white]    (-0.25,-0.25) rectangle (0.25,0.25);
  \fill[blue!18]  (-0.25,-0.25) rectangle (0.25,0.25);
  \fill[white]    (-0.125,-0.125) rectangle (0.125,0.125);
  \fill[orange!25](-0.125,-0.125) rectangle (0.125,0.125);
  \fill[white]    (-0.0625,-0.0625) rectangle (0.0625,0.0625);
  \foreach \r in {1,0.5,0.25,0.125,0.0625}{\draw[thin] (-\r,-\r) rectangle (\r,\r);}
  \draw[->] (-1.15,0) -- (1.2,0) node[right] {$x$};
  \draw[->] (0,-1.15) -- (0,1.2) node[above] {$y$};
  \draw[dashed,gray] (-1,-1) -- (1,1) node[above right,black] {\small $x=y$};
  \node at (-0.75,0.82) {\small $S_0$};
  \node at (-0.375,0.39) {\small $S_1$};
  \node at (-0.19,0.19) {\tiny $S_2$};
  \fill[red!70] (0.75,0.375) circle (0.05);
  \node[right] at (0.8,0.375) {\small $E$};
  \fill[red!70] (0.375,0.1875) circle (0.025);
  \node[right] at (0.40,0.17) {\scriptsize $\tfrac12E$};
  \fill[red!70] (0.1875,0.09375) circle (0.0125);
  \fill[red!70] (0.09375,0.046875) circle (0.00625);
  \node[below right] at (1,0) {\scriptsize $1$};
  \node[below right] at (0.5,0) {\scriptsize $\tfrac12$};
\end{tikzpicture}
\hspace{1.2cm}
\begin{tikzpicture}[scale=2.6]
  \fill[blue!18]   (0,0) rectangle (1,1);
  \fill[orange!25] (0,0) rectangle (0.5,0.5);
  \fill[blue!18]   (0,0) rectangle (0.25,0.25);
  \fill[orange!25] (0,0) rectangle (0.125,0.125);
  \fill[white]     (0,0) rectangle (0.0625,0.0625);
  \foreach \r in {1,0.5,0.25,0.125,0.0625}{\draw[thin] (\r,0) -- (\r,\r) -- (0,\r);}
  \draw[->] (0,0) -- (1.2,0) node[right] {$|x|$};
  \draw[->] (0,0) -- (0,1.2) node[above] {$|y|$};
  \draw[dashed,gray] (0,0) -- (1,1);
  \node at (0.75,0.85) {\small $S_0$};
  \node at (0.37,0.43) {\small $S_1$};
  \node at (0.19,0.215) {\scriptsize $S_2$};
  \draw[<-] (0.12,0.75) -- (-0.25,0.75) node[left,align=right] {\scriptsize $|x|$ small,\\[-2pt] \scriptsize $|y|\sim1$};
  \draw[<-] (0.75,0.12) -- (0.75,-0.22) node[below,align=center] {\scriptsize $|x|\sim1$, $|y|$ small};
  \node[below] at (1,0) {\scriptsize $1$};
  \node[below] at (0.5,0) {\scriptsize $\tfrac12$};
  \node[left] at (0,1) {\scriptsize $1$};
  \node[left] at (0,0.5) {\scriptsize $\tfrac12$};
\end{tikzpicture}
\caption{The dyadic shells $S_k=\{2^{-k-1}<\max\{|x|,|y|\}\le2^{-k}\}$. \emph{Left:} $d=1$, in the $(x,y)$-plane, where $S_k$ is the square frame between the squares of half-side $2^{-k-1}$ and $2^{-k}$ (shells alternately shaded); the small discs are the set $E$ of Step~4 around $(x_0,y_0)=(\frac34,\frac38)$ and its copies $2^{-k}E\subset S_k$. \emph{Right:} general $d$, drawn in the $(|x|,|y|)$-quadrant, where $S_k$ is an L-shaped band. Each shell contains three kinds of pairs: both points of size $\sim2^{-k}$ (near the corner, including the diagonal $x=y$, where the kernel is singular), and the two mixed parts where one point has size $\sim2^{-k}$ and the other is much closer to the origin (near the axes, where $|x|^a$ or $|y|^c$ is singular). All $S_k$ are copies of $S_0$ scaled by $2^{-k}$; only the joint singularity at $(0,0)$ is approached as $k\to\infty$.}
\label{fig:shells}
\end{figure}
 If $G\ge0$ is homogeneous of degree $\theta$ on $\R^{2d}$, i.e.\ $G(\lambda x,\lambda y)=\lambda^\theta G(x,y)$, then the substitution $(x,y)=2^{-k}(x',y')$ gives
\begin{equation}\label{eq:scaling}
  \int_{S_k} \int_{S_k} G\,dx\,dy=2^{-k(\theta+2d)}\int_{S_0} \int_{S_0} G\,dx'\,dy' .
\end{equation}

\medskip
\emph{Step 3: upper bound; $\gamma>0$ implies $\int_Q \int_Q F_u<\infty$.}  \ Notice that $m(u)=\infty$ if all derivatives of $u$ of order $\ge1$ vanish at $0$ (e.g.\ if $u$ is constant near $0$). If $m=m(u)<\infty$, then $P(x):=\sum_{|\beta|=m}\frac{\partial^\beta u(0)}{\beta!}x^\beta$ is a nonzero homogeneous polynomial of degree $m$, and Taylor's formula gives, for $|x|\le1$,
\begin{equation}\label{eq:taylor}
  u(x)=u(0)+P(x)+R(x),\qquad |R(x)|\le C_R|x|^{m+1},\qquad |\nabla R(x)|\le C_R|x|^{m}.
\end{equation}

Assume first $m<\infty$. By \eqref{eq:taylor}, for $|z|\le1$,
\[
  |\nabla u(z)|=|\nabla P(z)+\nabla R(z)|\le C|z|^{m-1}+C_R|z|^m\le C_1|z|^{m-1},
\]
since $\nabla P$ is homogeneous of degree $m-1$. For $x,y\in B_1$ every point $z$ of the segment $[x,y]$ satisfies $|z|\le\rho(x,y)$, hence by the mean value inequality
\begin{equation}\label{eq:upper}
  |u(x)-u(y)|\le C_1\,\rho(x,y)^{m-1}\,|x-y|\qquad(x,y\in B_1).
\end{equation}
Therefore $F_u\le C_1^p\,G$ on $Q$, where
\[
  G(x,y):=|x|^{a}|y|^{c}\,\rho(x,y)^{(m-1)p}\,|x-y|^{p-d-sp}
\]
is homogeneous of degree $\theta=a+c+(m-1)p+p-d-sp=\gamma-2d$. By \eqref{eq:scaling},
\[
  \int_Q \int_Q F_u\le C_1^p\sum_{k\ge0}2^{-k\gamma}\,J_0,\qquad J_0:=\int_{S_0} \int_{S_0} G .
\]
We check $J_0<\infty$. On $S_0$ we have $\rho\le1$, so $\rho^{(m-1)p}\le1$. If $|x|\le\frac14$, then $|y|>\frac12$ (since $\rho>\frac12$) and $|x-y|\ge\frac14$, so $G\le C|x|^a$, which is integrable over $B_{1/4}$ because $a>-d$; the case $|y|\le\frac14$ is symmetric (using $c>-d$). If $|x|,|y|\in(\frac14,1]$, then $|x|^a|y|^c\le C$ and $|x-y|^{p(1-s)-d}$ is integrable over $B_1\times B_1$ because $p(1-s)>0$. Hence $J_0<\infty$, and the geometric series converges since $\gamma>0$.

If $m=\infty$, then for every $N\in\N$ all derivatives of $u$ of order $\le N+1$ vanish at $0$, so $|\nabla u(z)|\le C_N|z|^N$ on $B_1$, and \eqref{eq:upper} holds with $m-1$ replaced by $N$. The argument above then applies with $\gamma$ replaced by $\gamma_N:=a+c+d+p(N+1-s)$, which is positive for $N$ large. Thus $\iint_QF_u<\infty$ whenever $m=\infty$.

\medskip
\emph{Step 4: lower bound; $\gamma\le0$ implies $\int_Q \int_Q F_u=\infty$.} Let $m<\infty$ (notice we do not have to consider $m=\infty$ here as $\gamma=\infty$ then). Since $P\not\equiv0$ is homogeneous, there is $x_0$ with $|x_0|=\frac34$ and $P(x_0)\ne0$. Put $y_0:=\frac12x_0$; then $P(y_0)=2^{-m}P(x_0)\ne P(x_0)$, $\rho(x_0,y_0)=\frac34$ and $|x_0-y_0|=\frac38$. By continuity there are $\eta>0$ and $c_0>0$ such that the set
\[
  E:=B_\eta(x_0)\times B_\eta(y_0)
\]
satisfies $E\subset S_0$, and for  $(x,y) \in E$: $|P(x)-P(y)|\ge2c_0$, $\ \frac14\le|x|,|y|\le1$, $\ |x-y|\le2$.

For $\lambda=2^{-k}$ and $(x,y)\in\lambda E\subset S_k$, write $(x,y)=\lambda(x',y')$ with $(x',y')\in E$. By homogeneity of $P$ and \eqref{eq:taylor},
\[
  |u(x)-u(y)|\ \ge\ \lambda^m|P(x')-P(y')|-|R(x)|-|R(y)|\ \ge\ 2c_0\lambda^m-2C_R\lambda^{m+1}\ \ge\ c_0\lambda^m
\]
as soon as $\lambda\le c_0/(2C_R)$, i.e.\ for $k\ge k_0$. Moreover on $\lambda E$
\[
  |x|^a|y|^c=\lambda^{a+c}|x'|^a|y'|^c\ge c_1\lambda^{a+c},\qquad |x-y|^{-d-sp}=\lambda^{-d-sp}|x'-y'|^{-d-sp}\ge2^{-d-sp}\lambda^{-d-sp},
\]
with $c_1:=\min\{4^{-a},1\}\min\{4^{-c},1\}>0$ (this covers both signs of $a,c$). Since $|\lambda E|=\lambda^{2d}|E|$,
\[
  \int_{\lambda E} \int_{\lambda E} F_u\ \ge\ c_1\,c_0^p\,2^{-d-sp}\,|E|\,\lambda^{a+c+mp-d-sp+2d}=c_2\,2^{-k\gamma}\qquad(k\ge k_0).
\]
The sets $\lambda E=2^{-k}E\subset S_k$ are pairwise disjoint, so
\[
  \int_Q \int_Q F_u\ \ge\ c_2\sum_{k\ge k_0}2^{-k\gamma}=\infty\qquad\text{if }\gamma\le0 .
\]
This finishes the proof of the lemma.
\end{proof}

\section{Failure of the P\'olya--Szeg\H{o} Inequality for Weighted Gagliardo Seminorms}\label{sec3}

In this section, we establish the validity and failure of the P\'olya--Szeg\H{o} inequality for weighted fractional Gagliardo seminorms in different regimes of the total weight $\alpha_1+\alpha_2$. We begin by recalling the notions of decreasing rearrangement and Schwarz symmetrization, along with the classical P\'olya--Szeg\H{o} inequality for the unweighted Gagliardo seminorm. We then establish the case $\alpha_1+\alpha_2<0$ and prove Theorem \ref{thm:main} by constructing a family of translated functions whose weighted seminorms tend to zero, while their Schwarz symmetrizations remain unchanged. Next, we consider the critical case $\alpha_1+\alpha_2=0$ and prove Theorem \ref{thm:zero}. Finally, we establish a counterexample in a particular setting with $\alpha_1+\alpha_2>0$ and prove Theorem \ref{thm:onesided}.

\smallskip

In the last part of this section, we establish a partial positive result in the case $p=1$, $\alpha_1=0$, and $\alpha_2>0$, and prove Proposition \ref{prop:partial-positive}. More precisely, using the coarea formula, we reduce the problem to a weighted fractional isoperimetric inequality and show that the P\'olya--Szeg\H{o} inequality holds for functions whose superlevel sets are balls, although their centres may vary with the level.

\smallskip

Let $u:\mathbb{R}^d\to\mathbb{R}$ be a measurable function which vanishes at infinity. The distribution function of $u$ is defined by
\begin{equation*}
  \lambda_u(t)
:=
\big|\{x\in\mathbb{R}^d:|u(x)|>t\}\big|,
\qquad \text{for } t>0.  
\end{equation*}
The decreasing rearrangement of $u$, denoted by $u_1^*$, is defined as
\begin{equation*}
   u_1^*(t)
:=
\inf\{s\geq0:\lambda_u(s)\leq t\},
\qquad \text{for } t\geq0. 
\end{equation*}
Using the decreasing rearrangement, the Schwarz symmetrization of $u$ is defined by
\begin{equation*}
   u^*(x)
:=
u_1^*\bigl(|B_1(0)|\,|x|^d\bigr),
\qquad \text{for } x\in\mathbb{R}^d. 
\end{equation*}

Thus, $u^*$ is a radially symmetric and radially decreasing function. Moreover, its superlevel sets are balls centred at the origin and have the same measure as the corresponding superlevel sets of $|u|$. In other words,
\begin{equation*}
    \big|\{x\in\mathbb{R}^d:u^*(x)>t\}\big|
=
\big|\{x\in\mathbb{R}^d:|u(x)|>t\}\big|,
\qquad \text{for } t>0.
\end{equation*}

Therefore, the Schwarz symmetrization preserves the distribution function of $u$. As a consequence, it also preserves all $L^q$-norms for which they are finite. In particular,
\begin{equation*}
    \|u^*\|_{L^q(\mathbb{R}^d)}
=
\|u\|_{L^q(\mathbb{R}^d)},
\qquad  1\leq q\leq\infty.
\end{equation*}
The importance of this rearrangement in the study of variational problems
comes from the fact that, for a large class of unweighted energies,
symmetrization does not increase the corresponding energy.

Recall that the fractional version of P\'olya--Szeg\H{o} inequality (see \cite[Section $9.2$]{Almgren1989} and \cite[Theorem A.$1.$, Appendix A]{frank2008non}),
\begin{equation*}
\int_{\mathbb{R}^d}\int_{\mathbb{R}^d}
\frac{|u^{*}(x)-u^{*}(y)|^p}
{|x-y|^{d+sp}}
\,dx\,dy
\leq
\int_{\mathbb{R}^d}\int_{\mathbb{R}^d}
\frac{|u(x)-u(y)|^p}
{|x-y|^{d+sp}}
\,dx\,dy.
\end{equation*}
Thus, in the unweighted setting, passing to the Schwarz symmetrization preserves the relevant $L^p$-information while decreasing the Gagliardo seminorm. This makes radially symmetric nonincreasing functions particularly useful in the study of fractional variational inequalities. However, the situation is different for the weighted Gagliardo seminorm defined in \eqref{Defn: weighted Gagliardo seminorm}. Unlike the unweighted seminorm, the weighted seminorm is not invariant under translations, since the weights are centred at the origin. This naturally raises the question of whether Schwarz symmetrization can still decrease the weighted seminorm, possibly up to a multiplicative constant.

\smallskip

\subsection{Failure of the P\'olya--Szeg\H{o} inequality for negative total weight} In this subsection, we prove Theorem \ref{thm:main} by constructing a suitable family of translations of a fixed radially symmetric function. The key observation is that, when $\alpha_1+\alpha_2<0$, the corresponding weighted Gagliardo seminorm tends to zero as the translation parameter tends to infinity, whereas the Schwarz symmetrization remains unchanged.

\begin{proof}[\textbf{Proof of Theorem \ref{thm:main}}]
\emph{Step 1: $u_R^*=\varphi$, and $\sn{\varphi}$ is finite and positive.} Since translation does not change the distribution function and $\varphi$ is radially symmetric and nonincreasing, its Schwarz symmetrization remains unchanged. Hence,
\begin{equation*}
    u_R^*=\varphi^*=\varphi.
\end{equation*}
As $\varphi \in C^{\infty}_{c}(\mathbb{R}^{d})$, its weighted Gagliardo seminorm is finite. Indeed, by Lemma \ref{lemma:critical interval}, we have
\begin{equation*}
    [\varphi]_{W^{s,p,\alpha}(\mathbb{R}^{d})}<\infty.
\end{equation*}

\medskip
\emph{Step 2: change of variables.}  Let $x_0=Re$ and make the change of variables $x=x_0+\xi$ and
$y=x_0+\eta$ in the weighted Gagliardo seminorm of $u_R$. We then obtain
\begin{equation*}
[u_R]^{p}_{W^{s,p, \alpha}(\mathbb{R}^{d})}=\int_{\mathbb{R}^{d}}\int_{\mathbb{R}^{d}}\frac{|x_0+\xi|^{\alpha_{1}p}\,|x_0+\eta|^{\alpha_{2}p}\,|\varphi(\xi)-\varphi(\eta)|^p}{|\xi-\eta|^{d+sp}}\,d\xi\,d\eta .
\end{equation*}
As $\operatorname{supp}\varphi\subset B_2$ and $R\geq 8$, the support of
the integrand can be covered by the following three regions:
\begin{equation*}
    \mathrm{I}=\{|\xi|,|\eta|<R/2\},\qquad
  \mathrm{II}=\{\xi\in B_2,\ |\eta|\ge R/2\},\qquad
  \mathrm{III}=\{\eta\in B_2,\ |\xi|\ge R/2\}. 
\end{equation*}
We denote the corresponding parts of the integral by the same letters.

\smallskip
\medskip
\emph{Step 3: region I (near field).}
For $|\xi|<R/2$, we have $R/2\le|x_0+\xi|\le 3R/2$, hence
\begin{equation*}
   |x_0+\xi|^{\alpha_{1} p}\le c_{\alpha_{1}} R^{\alpha_{1}p} 
\end{equation*} 
with $c_{\alpha_{1}}=\max\{2^{-\alpha_{1}p},(3/2)^{\alpha_{1}p}\}$, and similarly, $|x_0+\eta|^{\alpha_{2}p}\leq c_{\alpha_{2}} R^{\alpha_{2}p}$. Therefore,
\begin{equation*}
     \mathrm{I}\le c_{\alpha_{1}} c_{\alpha_{2}} \,R^{\alpha_{1}p+\alpha_{2}p}\,[\varphi]^p_{W^{s,p}(\mathbb{R}^{d})},
\end{equation*}
where $[\varphi]_{W^{s,p}(\mathbb{R}^{d})}<\infty$ is the unweighted Gagliardo seminorm of $\varphi\in C_c^\infty(\mathbb{R}^{d})$.

\smallskip

\medskip
\emph{Step 4: region II (cross term).}
On $\mathrm{II}$, since $|\eta|\geq R/2\geq 4$, we have
\begin{equation*}
    |\xi-\eta|\geq |\eta|-|\xi|\geq |\eta|-2\geq \frac{|\eta|}{2}.
\end{equation*}
Moreover, $ |\varphi(\xi)-\varphi(\eta)|\leq 1$ and
\begin{equation*}
    |x_0+\xi|^{\alpha_1p}
\leq c_{\alpha_1}R^{\alpha_1 p}.
\end{equation*}
Therefore,
\begin{equation*}
\mathrm{II}
\leq
c_{\alpha_1}\,2^{d+sp}\,|B_2|\,R^{\alpha_1p}
\int_{|y-x_0|\geq R/2}
\frac{|y|^{\alpha_2p}}{|y-x_0|^{d+sp}}\,dy,
\end{equation*}
where we have written $y=x_0+\eta$. Since $x_0=Re$, we have
$|\eta|=|y-x_0|$. Rescaling by $y=Rz$, we obtain
\begin{equation*}
\mathrm{II}
\leq
C R^{\alpha_1p}R^{\alpha_2p-sp}K
=
C R^{\alpha_1p+\alpha_2p-sp}K,
\end{equation*}
where
\begin{equation*}
K:=
\int_{|z-e|\geq 1/2}
\frac{|z|^{\alpha_2p}}{|z-e|^{d+sp}}\,dz.
\end{equation*}
The integral $K$ is finite. Indeed, near $z=0$ this follows from
$\alpha_2p>-d$, while at infinity the integrand behaves like
$|z|^{\alpha_2p-d-sp}$, which is integrable since $\alpha_2p<sp$.
Moreover, the possible singularity at $z=e$ is excluded by the condition
$|z-e|\geq 1/2$. Hence,
\begin{equation*}
\mathrm{II}
\leq
C R^{\alpha_1p+\alpha_2p-sp}.
\end{equation*}

\smallskip

\medskip
\emph{Step 5: region III.}
By a similar argument, interchanging the roles of $\alpha_{1}p$ and $\alpha_{2}p$ and using
the assumptions $\alpha_{1}p>-d$ and $\alpha_{1} p <sp$, we obtain
\begin{equation*}
    \mathrm{III}\leq C\,R^{\alpha_{1}p+\alpha_{2}p-sp}.
\end{equation*}

\medskip
\emph{Step 6: conclusion.}
Adding the above three inequalities, we obtain
\begin{equation*}
    [u_R]^{p}_{W^{s,p, \alpha}(\mathbb{R}^{d})}\le C\left(R^{(\alpha_1+\alpha_2)p}+R^{(\alpha_1+\alpha_2-s)p}\right)\longrightarrow 0\quad \text{as }  R\to\infty ,
\end{equation*}
since $\alpha_1+\alpha_2<0$. Suppose, on the other hand, that there exists a constant $C>0$ such that $[u^*]_{W^{s,p,\alpha}(\mathbb{R}^{d})}
\leq
C[u]_{W^{s,p,\alpha}(\mathbb{R}^{d})}$ holds for every nonnegative $u\in C_c^\infty(\mathbb{R}^{d})$. Then, using
$u_R^*=\varphi$, we would have
\begin{equation*}
  0<[\varphi]_{W^{s,p, \alpha}(\mathbb{R}^{d})}=[u_R^*]_{W^{s,p, \alpha}(\mathbb{R}^{d})}\le C\,[u_R]_{W^{s,p, \alpha}(\mathbb{R}^{d})}\longrightarrow 0,  
\end{equation*}
which is a contradiction. Hence, no such constant $C$ can exist.
\end{proof}

\subsection{The case of zero total weight: failure with  \texorpdfstring{$\alpha_1+\alpha_2 =0$}{alpha1+alpha2=0}}

In this subsection, we prove Theorem \ref{thm:zero}, which establishes the failure of the Pólya–Szegő inequality in the case of zero total weight, $\alpha_1+\alpha_2=0$. The proof relies on an exact symmetrized representation of the weighted Gagliardo seminorm and a family of translations of a fixed radially symmetric function. We show that the seminorm of the translated function converges to the corresponding unweighted fractional Gagliardo seminorm, whereas the seminorm of its Schwarz symmetrization exceeds this limit by a strictly positive quantity.

When $\alpha_1+\alpha_2=0$, the weight is homogeneous of degree $0$, and the seminorm admits an exact symmetrized form.

\medskip
\noindent\textbf{Exact formula.}
Let $\alpha_1=\beta$, $\alpha_2=-\beta$. Since $|u(x)-u(y)|^p|x-y|^{-d-sp}$ is symmetric in $(x,y)$, we may replace $|x|^{\beta p}|y|^{-\beta p}$ by $\tfrac12\big(t+t^{-1}\big)=\cosh(\log t)$ with $t=(|x|/|y|)^{\beta p}$. Hence, for every measurable $u$,
\begin{equation}\label{eq:cosh}
  \sn{u}^p=\int_{\R^d}\int_{\R^d}\cosh\!\Big(\beta p\,\log\frac{|x|}{|y|}\Big)\,\frac{|u(x)-u(y)|^p}{|x-y|^{d+sp}}\,dx\,dy .
\end{equation}
In particular $\sn{u}\ge[u]_{W^{s,p}(\R^d)}$, since $\cosh\ge1$ with equality exactly on $\{|x|=|y|\}$. For $u_R(x)=\varphi(x-Re)$ and $x_0=Re$, \eqref{eq:cosh} gives
\[
  \sn{u_R}^p=\int_{\R^d}\int_{\R^d}\cosh\!\Big(\beta p\,\log\frac{|x_0+\xi|}{|x_0+\eta|}\Big)\frac{|\varphi(\xi)-\varphi(\eta)|^p}{|\xi-\eta|^{d+sp}}\,d\xi\,d\eta ,
\]
while for $u_R^*=\varphi$,
\begin{align*}
    \sn{\varphi}^p  =[\varphi]^p_{W^{s,p}(\R^d)}+\delta_\beta(\varphi),
\end{align*}
where 
\begin{equation}\label{Defn: delta_beta}
    \ 
  \delta_\beta(\varphi) :=\int_{\R^d}\int_{\R^d}\Big(\cosh\!\Big(\beta p\log\frac{|x|}{|y|}\Big)-1\Big)\frac{|\varphi(x)-\varphi(y)|^p}{|x-y|^{d+sp}}\,dx\,dy .
\end{equation}

\bigskip

We now prove Theorem \ref{thm:zero} by comparing the weighted Gagliardo seminorms of a translated function and its Schwarz symmetrization.

\begin{proof}[\textbf{Proof of Theorem \ref{thm:zero}}]
The condition $|\beta|<\min\{s,d/p\}$ is exactly $-d/p<\alpha_1,\alpha_2<s$, so by Step~1 of the proof of Theorem~\ref{thm:main} all quantities involved are finite.

\medskip
\emph{(1) Positivity of $\delta_\beta(\varphi)$.}
Since $\varphi$ is radial, $\varphi(x)\neq\varphi(y)$ forces $|x|\neq|y|$, and there $$\cosh(\beta p\log(|x|/|y|))>1$$because $\beta\neq0$. The set $\{(x,y):\varphi(x)\neq\varphi(y)\}$ contains $B_1\times B_2^c$, which has positive measure. Hence the integrand of $\delta_\beta(\varphi)$ is positive on a set of positive measure, and $\delta_\beta(\varphi)>0$.

\medskip
\emph{(2) The limit.}
Use the regions $\mathrm I,\mathrm{II},\mathrm{III}$ of Step~2 in the proof of Theorem~\ref{thm:main}, with $R\ge 8$.

On $\mathrm I=\{|\xi|,|\eta|<R/2\}$, both $|x_0+\xi|$ and $|x_0+\eta|$ lie in $[R/2,3R/2]$, so their ratio lies in $[1/3,3]$ and the $\cosh$ factor is bounded by $\cosh(\beta p\log3)$. This gives the $R$-independent integrable majorant
\[
  \cosh(\beta p\log 3)\,\frac{|\varphi(\xi)-\varphi(\eta)|^p}{|\xi-\eta|^{d+sp}}\in L^1(\R^d\times\R^d).
\]
For fixed $(\xi,\eta)$, we have $|x_0+\xi|/|x_0+\eta|\to1$, so the $\cosh$ factor tends to $1$; moreover $\mathbf 1_{\mathrm I}\to1$ pointwise. By dominated convergence, the contribution of $\mathrm I$ tends to $[\varphi]^p_{W^{s,p}(\R^d)}$.

On $\mathrm{II}$ and $\mathrm{III}$, Steps~4 and~5 of the proof of Theorem~\ref{thm:main} with $(\alpha_{1}+\alpha_{2})p=0$ give the bound $C\,R^{\alpha_{1}p+\alpha_{2}p-sp}=C\,R^{-sp}\to0$.

\medskip
\emph{(3) Conclusion.}
Since $u_R^*=\varphi$, parts (1) and (2) give
\begin{align*}
    \sn{u_R^*}^p-\sn{u_R}^p & =\sn{\varphi}^p-\sn{u_R}^p \\ & \quad \longrightarrow \sn{\varphi}^p-[\varphi]^p_{W^{s,p}(\R^d)}=\delta_\beta(\varphi)>0 .
\end{align*}
This finishes the proof of the theorem.\end{proof}



\subsection{Failure of the P\'olya--Szeg\H{o} Inequality for \texorpdfstring{$\alpha_1=0$, $\alpha_2=\varepsilon>0$}{alpha1=0, alpha2=epsilon>0}}

In this subsection, we prove Theorem~\ref{thm:onesided}, which establishes the failure of the P\'olya--Szeg\H{o} inequality when $\alpha_1=0$ and $\alpha_2=\varepsilon>0$. Let $\alpha_1=0$ and $\alpha_2=\varepsilon>0$, so that $\alpha=\varepsilon>0$, and write $b=\varepsilon p$. The corresponding seminorm is given by
\[
[u]_{0,\varepsilon}^p=\int_{\R^d}\int_{\R^d}\frac{|y|^{b}\,|u(x)-u(y)|^p}{|x-y|^{d+sp}}\,dx\,dy .
\]
It is finite for $\varphi\in C_c^\infty(\R^d)$ when $\varepsilon<s$, by Step~1 of the proof of Theorem~\ref{thm:main}. Recall (as in introduction section),

\[
  m_\varphi(\eta):=\int_{\R^d}\frac{|\varphi(\xi)-\varphi(\eta)|^p}{|\xi-\eta|^{d+sp}}\,d\xi .
\]
It is bounded, nonnegative, and radial if $\varphi$ is radial. Substituting $x=\xi+h$, $y=\eta+h$ gives, for every $h\in\R^d$,
\begin{equation}\label{eq:marg}
  [\varphi(\cdot-h)]_{0,\varepsilon}^p=\int_{\R^d}|\eta+h|^{b}\,m_\varphi(\eta)\,d\eta .
\end{equation}


\bigskip

\begin{proof}[\textbf{Proof of Theorem \ref{thm:onesided}}]
Write $m=m_\varphi$ and $\tilde{\Delta}(h)=[\varphi(\cdot-h)]_{0,\varepsilon}^p-[\varphi]_{0,\varepsilon}^p$.

\medskip
\emph{Step 1: properties of $m$.} The marginal $m$ is even and bounded, since near the diagonal $|\varphi(\xi)-\varphi(\eta)|^p\le L^p|\xi-\eta|^p$ and $p-1-sp>-1$. For $|\eta|<1/2$ we have $\varphi(\eta)=1$, so
\[
  m(\eta)=\int_{|\xi|>1}|\varphi(\xi)-1|^p\,|\xi-\eta|^{-1-sp}\,d\xi .
\]
This is a smooth even function of $\eta$ on $(-\tfrac12,\tfrac12)$, with $|m''|\le C_1:=(1+sp)(2+sp)\int_{|\xi|>1}(|\xi|-\tfrac12)^{-3-sp}d\xi$. Hence $|m(\eta)-m(0)|\le C_1\eta^2$ there. Note that $C_1$ depends only on $s,p$ and not on $\varphi$. Together with boundedness of $m$, this shows that $Q_b(\varphi)$ converges absolutely.

\medskip
\emph{Step 2: reduction.} By \eqref{eq:marg} and evenness of $m$,
\[
  \Delta(h)=\int_{\R}D_h(\eta)\,m(\eta)\,d\eta,\qquad D_h(\eta):=\tfrac12\big(|\eta+h|^b+|\eta-h|^b\big)-|\eta|^b .
\]
For $0<b<1$, $D_h$ is integrable and $\int_\R D_h=0$. Indeed, with $F(t)=t^{b+1}/(b+1)$, a direct computation gives
\[
  \int_{-A}^{A}D_h=F(A+h)-2F(A)+F(A-h)\le h^2\,b\,(A-h)^{b-1}\to0\quad(A\to\infty).
\]
Therefore
\[
  \Delta(h)=\int_\R D_h(\eta)\big(m(\eta)-m(0)\big)\,d\eta .
\]

\medskip
\emph{Step 3: asymptotics.} Let $0<h<1/4$.
\begin{itemize}
\item On $\{|\eta|<2h\}$ we have $|\eta\pm h|<3h$ and $|\eta|<2h$. Since both terms of $D_h$ are nonnegative,
\[
|D_h(\eta)|\le \max\Big\{\tfrac12\big(|\eta+h|^b+|\eta-h|^b\big),\,|\eta|^b\Big\}\le (3h)^b .
\]
(The sign of $D_h$ is not fixed on this set; only this bound is used.) Together with
$|m(\eta)-m(0)|\le C_1\eta^2$ for $|\eta|<\tfrac12$ (valid since $2h<\tfrac12$), the contribution of this region is at most
\[
(3h)^b\,C_1\int_{|\eta|<2h}\eta^2\,d\eta=\tfrac{16}{3}\,3^b C_1\,h^{3+b}=O(h^{3+b})=o(h^2).
\]
  \item On $\{|\eta|\ge2h\}$, Taylor's formula applied to $t\mapsto t^b$ on $[|\eta|-h,|\eta|+h]\subset(0,\infty)$ gives
  \[
    |D_h(\eta)|\le\tfrac{h^2}{2}\,b(1-b)\,(|\eta|-h)^{b-2}\le C h^2|\eta|^{b-2},
    \qquad
    \frac{D_h(\eta)}{h^2}\to\frac{b(b-1)}{2}|\eta|^{b-2}.
  \]
  The function $|\eta|^{b-2}|m(\eta)-m(0)|$ is integrable by Step~1, so dominated convergence applies.
\end{itemize}
This proves the stated limit. If $Q_b(\varphi)>0$ then, since $b(b-1)<0$, we get $\Delta(h)<0$ for all small $h>0$. Also $\varphi(\cdot-h)^*=\varphi$.

\medskip
\emph{Step 4: a profile with $Q_b>0$ for all $b\in(0,\tfrac12]$.} For $\delta\in(0,1]$, let $\varphi_\delta\in C_c^\infty(\R)$ be even, nonincreasing on $[0,\infty)$, equal to $1$ on $[-1,1]$, and supported in $[-1-\delta,1+\delta]$, with $\varphi_\delta\to\mathbf 1_{[-1,1]}$ pointwise as $\delta\to0$. Write $m_\delta=m_{\varphi_\delta}$ and split $Q_b(\varphi_\delta)$ into three ranges.
\begin{itemize}
  \item $|\eta|<\tfrac12$: by Step~1 the contribution is at least $-C_1\int_{|\eta|<1/2}|\eta|^b\,d\eta\ge-C_1$.
  \item $|\eta|>2$: here $m_\delta-m_\delta(0)\ge-m_\delta(0)$, and $m_\delta(0)\le\int_{|\xi|>1}|\xi|^{-1-sp}\,d\xi=\tfrac2{sp}$. The contribution is at least $-\tfrac2{sp}\int_{|\eta|>2}|\eta|^{b-2}\,d\eta\ge-\tfrac{8}{sp}$ for $b\le\tfrac12$.
  \item $\tfrac12\le|\eta|\le2$: here $|\eta|^{b-2}\ge 2^{b-2}\ge\tfrac14$, and, for $b\in(0,\tfrac12]$,
$\int_{\{1/2\le|\eta|\le2\}}|\eta|^{b-2}\,d\eta\le 3$. Since $m_\delta(0)\le \frac{2}{sp}$,
the contribution is at least
\[
\frac14\int_{1/2\le|\eta|\le2} m_\delta\,d\eta \;-\; 3\,m_\delta(0)
\;\ge\; \frac14\int_J\!\int_J \frac{|\varphi_\delta(\xi)-\varphi_\delta(\eta)|^p}{|\xi-\eta|^{1+sp}}\,d\xi\,d\eta \;-\; \frac{6}{sp},
\quad J=(\tfrac12,2).
\]
\end{itemize}
By Fatou's lemma, as $\delta\to0$ the double integral over $J\times J$ has $\liminf$ at least
\[
  2\int_{1/2}^{1}\int_{1}^{2}|\xi-\eta|^{-1-sp}\,d\xi\,d\eta=\infty,\qquad\text{since } sp\ge1 .
\]
Hence $Q_b(\varphi_\delta)\to\infty$ as $\delta\to0$, uniformly in $b\in(0,\tfrac12]$. Fixing $\delta$ small and setting $\varphi=\varphi_\delta$ gives $Q_b(\varphi)>0$ for all $b=\varepsilon p\in(0,\tfrac12]$. Step~3 then applies for every $\varepsilon\in(0,\varepsilon_0]$, where $\varepsilon<s$ ensures finiteness.
\end{proof}

\subsection{Proof of the Partial result for \texorpdfstring{$p=1$}{sp=1}}

In this subsection, we establish a partial P\'olya--Szeg\H{o} inequality for $p=1$. We first derive a coarea formula that reduces the problem to a weighted fractional isoperimetric inequality. We then prove the desired inequality for functions whose superlevel sets are balls.

\smallskip

The next lemma establishes a coarea formula for the case $p=1$, expressing the weighted fractional seminorm of a nonnegative function in terms of the weighted perimeters of its superlevel sets. Let $P_b$ be as defined in \eqref{Defn: Pb}. 

\begin{lemma}[Coarea formula]\label{lem:coarea}
Let $p=1$ and let $u\ge0$ be measurable. Then
\begin{equation}\label{eq:coarea}
  [u]_{0,\varepsilon}=\int_0^\infty P_b(\{u>t\})\,dt\qquad\text{in }[0,\infty].
\end{equation}
\end{lemma}

\begin{proof}
\emph{Step 1: a pointwise layer-cake identity.} For real numbers $\alpha,\beta\ge0$,
\begin{equation}\label{eq:layer}
  |\alpha-\beta|=\int_0^\infty\big|\mathbf 1_{(t,\infty)}(\alpha)-\mathbf 1_{(t,\infty)}(\beta)\big|\,dt .
\end{equation}
Indeed, if $\alpha\ge\beta$, then $\mathbf 1_{(t,\infty)}(\alpha)-\mathbf 1_{(t,\infty)}(\beta)$ equals $1$ when $\beta\le t<\alpha$ and $0$ for all other $t\ge0$, so the integral is the length $\alpha-\beta$ of the interval $[\beta,\alpha)$. The case $\alpha<\beta$ is symmetric. Applying \eqref{eq:layer} with $\alpha=u(x)$ and $\beta=u(y)$, and noting $\mathbf 1_{(t,\infty)}(u(x))=\mathbf 1_{\{u>t\}}(x)$, gives
\[
  |u(x)-u(y)|=\int_0^\infty\big|\mathbf 1_{\{u>t\}}(x)-\mathbf 1_{\{u>t\}}(y)\big|\,dt\qquad\text{for all }x,y\in\R^d .
\]
Geometrically, the pair $(x,y)$ is separated by the level set $\{u>t\}$ (one point inside, one outside) exactly for $t$ between $u(y)$ and $u(x)$.

\emph{Step 2: measurability.} The set $\{(x,t)\in\R^d\times(0,\infty):u(x)>t\}$ is measurable because $u$ is measurable. Hence $(x,y,t)\mapsto|\mathbf 1_{\{u>t\}}(x)-\mathbf 1_{\{u>t\}}(y)|\,|y|^b|x-y|^{-d-s}$ is a nonnegative measurable function on $\R^d\times\R^d\times(0,\infty)$.

\emph{Step 3: Tonelli.} Multiply the identity of Step~1 by $|y|^b|x-y|^{-d-s}$, integrate over $(x,y)\in\R^d\times\R^d$, and exchange the order of integration (allowed by Tonelli's theorem since everything is nonnegative):
\[
  [u]_{0,\varepsilon}=\int_0^\infty\Big(\int_{\R^d}\int_{\R^d}\frac{|y|^b\,|\mathbf 1_{\{u>t\}}(x)-\mathbf 1_{\{u>t\}}(y)|}{|x-y|^{d+s}}\,dx\,dy\Big)\,dt=\int_0^\infty P_b(\{u>t\})\,dt .
\]
Here the inner integral is $[\mathbf 1_{\{u>t\}}]_{0,\varepsilon}=P_b(\{u>t\})$ by definition. If $|\{u>t\}|=\infty$ for some $t$, the inner integral is still defined in $[0,\infty]$, and the identity remains valid.
\end{proof}

Note that \eqref{eq:coarea} is special to $p=1$. For $p>1$ the map $t\mapsto|\mathbf 1_{\{u>t\}}(x)-\mathbf 1_{\{u>t\}}(y)|$ still integrates to $|u(x)-u(y)|$, but not to $|u(x)-u(y)|^p$, so the seminorm no longer splits into contributions of individual level sets.

\begin{lemma}[Reduction to sets]\label{cor:sets}
Let $p=1$. The following are equivalent:
\begin{enumerate}
  \item[(i)] $[u^*]_{0,\varepsilon}\le[u]_{0,\varepsilon}$ for every measurable $u$ vanishing at infinity (i.e.\ $|\{|u|>t\}|<\infty$ for all $t>0$), where $u^*:=|u|^*$;
  \item[(ii)] the weighted fractional isoperimetric inequality
  \[
    P_b(E)\ \ge\ P_b(E^*)\qquad\text{for every measurable } E\subset\R^d \text{ with } |E|<\infty ,
  \]
  where $E^*$ is the open ball centred at $0$ with $|E^*|=|E|$.
\end{enumerate}
\end{lemma}

\begin{proof}
\emph{(i)$\Rightarrow$(ii).} Take $u=\mathbf 1_E$. Then $u^*=\mathbf 1_{E^*}$, and by definition $[u]_{0,\varepsilon}=P_b(E)$ and $[u^*]_{0,\varepsilon}=P_b(E^*)$.

\emph{(ii)$\Rightarrow$(i).}
First we reduce to $u\ge0$. Since $\big||u(x)|-|u(y)|\big|\le|u(x)-u(y)|$, we have $[|u|]_{0,\varepsilon}\le[u]_{0,\varepsilon}$, and $u^*=|u|^*$ by definition. So it suffices to prove $[v^*]_{0,\varepsilon}\le[v]_{0,\varepsilon}$ for $v=|u|\ge0$.

Next we use a basic property of the symmetric decreasing rearrangement: its superlevel sets are the rearranged superlevel sets,
\[
  \{v^*>t\}=\{v>t\}^*\qquad\text{for every }t>0 .
\]
This follows from the layer-cake definition $v^*(x)=\int_0^\infty\mathbf 1_{\{v>\tau\}^*}(x)\,d\tau$: the sets $\{v>\tau\}^*$ are centred open balls that shrink as $\tau$ increases, so $v^*(x)>t$ holds exactly when $x\in\{v>t\}^*$ (see e.g.\ Lieb--Loss, \emph{Analysis}, Chapter~3). In particular $|\{v^*>t\}|=|\{v>t\}|<\infty$.

Applying \eqref{eq:coarea} to $v^*$ and to $v$, and (ii) to $E=\{v>t\}$ for each $t>0$,
\[
  [v^*]_{0,\varepsilon}=\int_0^\infty P_b(\{v^*>t\})\,dt=\int_0^\infty P_b\big(\{v>t\}^*\big)\,dt\ \le\ \int_0^\infty P_b(\{v>t\})\,dt=[v]_{0,\varepsilon}. \qedhere
\]
\end{proof}



\begin{proof}[\textbf{Proof of the  partial result: Proposition \ref{prop:partial-positive}}]
The proof has two parts. First we prove the inequality for a single ball (Steps 1--4): moving a ball away from the origin never decreases $P_b$. Then we pass to general $u$ through the coarea formula (Steps 5--6). We may assume $[u]_{0,\varepsilon}<\infty$, since otherwise there is nothing to prove.

\medskip
\emph{Step 1: $P_b$ as the weight tested against a marginal.}
For a measurable set $E$ with $|E|<\infty$ define its energy marginal
\[
  m_E(y):=\int_{\R^d}\frac{|\mathbf 1_E(x)-\mathbf 1_E(y)|}{|x-y|^{d+s}}\,dx
  =\begin{cases}\displaystyle\int_{E^c}\frac{dx}{|x-y|^{d+s}}, & y\in E,\\[2ex] \displaystyle\int_{E}\frac{dx}{|x-y|^{d+s}}, & y\notin E.\end{cases}
\]
Thus $m_E(y)$ measures how strongly the point $y$ interacts with the side of $\partial E$ it does not belong to. By Tonelli's theorem, applied to the nonsymmetric form of the seminorm (the weight is $|y|^b$),
\begin{equation}\label{eq:PbE}
  P_b(E)=\int_{\R^d}\int_{\R^d}\frac{|y|^b\,|\mathbf 1_E(x)-\mathbf 1_E(y)|}{|x-y|^{d+s}}\,dx\,dy=\int_{\R^d}|y|^b\,m_E(y)\,dy .
\end{equation}
So $P_b(E)$ is the weight $|y|^b$ integrated against the density $m_E$.

\medskip
\emph{Step 2: the marginal of a ball, and translation.}
Let $E=B_\rho=B_\rho(0)$ and write $m:=m_{B_\rho}$. We record three properties of $m$.
\begin{enumerate}
  \item[(a)] $m$ is radial. If $R$ is a rotation, then $RB_\rho=B_\rho$ and $|Rx-Ry|=|x-y|$, so the substitution $x\mapsto Rx$ gives $m(Ry)=m(y)$.
  \item[(b)] $m\in L^1(\R^d)$, since $\int_{\R^d}m=2\int_{B_\rho}\int_{B_\rho^c}|x-y|^{-d-s}\,dx\,dy=2P_s(B_\rho)<\infty$ for $s\in(0,1)$.
  \item[(c)] $m(y)\le C\,|y|^{-d-s}$ for $|y|\ge2\rho$, because then $|x-y|\ge|y|/2$ for all $x\in B_\rho$, so $m(y)\le|B_\rho|\,2^{d+s}|y|^{-d-s}$.
\end{enumerate}
For a translated ball $B_\rho(h)=h+B_\rho$, the substitution $x\mapsto x+h$ in the definition of $m_E$ gives $m_{B_\rho(h)}(y)=m(y-h)$. Hence, by \eqref{eq:PbE} and the substitution $y=\eta+h$,
\begin{equation}\label{eq:Pbtrans}
  P_b(B_\rho(h))=\int_{\R^d}|y|^b\,m(y-h)\,dy=\int_{\R^d}|\eta+h|^b\,m(\eta)\,d\eta .
\end{equation}
In words: translating the ball by $h$ is the same as keeping the ball fixed and translating the weight by $-h$. The integral is finite: by (b) the region $|\eta|\le2\rho$ contributes at most $(2\rho+|h|)^b\|m\|_{L^1}$, and by (c) the region $|\eta|>2\rho$ contributes at most $C\int_{|\eta|>2\rho}(|\eta|+|h|)^b|\eta|^{-d-s}\,d\eta<\infty$, because $b<s$.

\medskip
\emph{Step 3: averaging over the sphere of centres.}
Put $f(h):=P_b(B_\rho(h))$. Then $f$ is radial. Indeed, for a rotation $R$, the substitution $\eta=R\eta'$ in \eqref{eq:Pbtrans}, together with (a), gives
\[
  f(Rh)=\int|R\eta'+Rh|^b\,m(R\eta')\,d\eta'=\int|\eta'+h|^b\,m(\eta')\,d\eta'=f(h).
\]
This is the only place where we use that the set is a ball: for a general set $E$ the marginal $m_E$ is not radial, and $P_b(E+h)$ need not be a radial function of $h$. Fix $r=|h|>0$ and let $\sigma_r$ be the normalized surface measure on $\{|h'|=r\}$. Since $f$ is constant on this sphere, Tonelli's theorem gives
\begin{equation}\label{eq:Pbavg}
  f(h)=\int f(h')\,d\sigma_r(h')=\int_{\R^d}M_r(\eta)\,m(\eta)\,d\eta,\qquad M_r(\eta):=\int|\eta+h'|^b\,d\sigma_r(h').
\end{equation}
Here $M_r(\eta)$ is the average of the weight $z\mapsto|z|^b$ over the sphere of radius $r$ centred at $\eta$. So moving the ball off the origin in a random direction replaces the weight $|\eta|^b$ by its spherical mean $M_r(\eta)$.

\medskip
\emph{Step 4: the spherical mean inequality $M_r(\eta)\ge|\eta|^b$.}
This is where $d\ge2$ is used. We show that $g(z)=|z|^b$ satisfies the sub-mean-value property. Since $g$ is not smooth at the origin, we regularize: for $\delta>0$ put $g_\delta(z):=(|z|^2+\delta^2)^{b/2}\in C^\infty(\R^d)$. A direct computation gives
\[
  \Delta g_\delta(z)=b\,(|z|^2+\delta^2)^{\frac b2-2}\Big((d+b-2)|z|^2+d\,\delta^2\Big)\ \ge\ 0\qquad(d\ge2,\ b>0),
\]
so $g_\delta$ is subharmonic. For a smooth function, the spherical mean $\mu(t):=\int g_\delta(\eta+th')\,d\sigma_1(h')$ satisfies, by the divergence theorem,
\[
  \mu'(t)=\int\nabla g_\delta(\eta+th')\cdot h'\,d\sigma_1(h')=\frac{1}{|\partial B_t|}\int_{B_t(\eta)}\Delta g_\delta\,dz\ \ge\ 0 ,
\]
and $\mu(t)\to g_\delta(\eta)$ as $t\to0^+$. Hence $\mu(r)\ge g_\delta(\eta)$, i.e.
\[
  \int\big(|\eta+h'|^2+\delta^2\big)^{b/2}\,d\sigma_r(h')\ \ge\ (|\eta|^2+\delta^2)^{b/2}.
\]
Letting $\delta\to0^+$ (the integrands decrease to $|\eta+h'|^b$ and are bounded on the sphere) yields
\begin{equation}\label{eq:submean}
  M_r(\eta)\ \ge\ |\eta|^b\qquad\text{for all }\eta\in\R^d,\ r>0 .
\end{equation}
Geometrically, $|z|^b$ grows from the origin in every direction, and in dimension $d\ge2$ this growth is strong enough (the term $(d-1)/|z|$ coming from the curvature of spheres) to beat the concavity of $t\mapsto t^b$ in the radial direction. In $d=1$ the curvature term is absent, $\Delta g=b(b-1)|z|^{b-2}<0$, and \eqref{eq:submean} fails; this is exactly the mechanism behind Theorem~\ref{thm:onesided}.

Inserting \eqref{eq:submean} into \eqref{eq:Pbavg} and using $m\ge0$ and \eqref{eq:PbE} for $E=B_\rho$, we obtain
\begin{equation}\label{eq:oneball}
  P_b(B_\rho(h))\ \ge\ \int_{\R^d}|\eta|^b\,m(\eta)\,d\eta=P_b(B_\rho)\qquad\text{for all }h\in\R^d,\ \rho>0 .
\end{equation}

\medskip
\emph{Step 5: the coarea formula.}
For $u\ge0$ and $x,y\in\R^d$,
\[
  |u(x)-u(y)|=\int_0^\infty\big|\mathbf 1_{\{u>t\}}(x)-\mathbf 1_{\{u>t\}}(y)\big|\,dt .
\]
Indeed, if, say, $u(x)\ge u(y)$, then the integrand equals $1$ exactly when $u(y)\le t<u(x)$ and $0$ otherwise. Multiplying by $|y|^b|x-y|^{-d-s}$ and applying Tonelli's theorem (everything is nonnegative) gives \eqref{eq:coarea}:
\[
  [u]_{0,\varepsilon}=\int_0^\infty P_b(\{u>t\})\,dt ,
\]
and the same identity holds for $u^*$. Note that $P_b(E)$ does not change if $E$ is modified on a null set.

\medskip
\emph{Step 6: comparing level sets.}
By hypothesis, for a.e.\ $t>0$ we have $\{u>t\}=B_{\rho(t)}(h(t))$ up to a null set. By definition of the symmetric decreasing rearrangement, $\{u^*>t\}$ is the centred ball with the same measure as $\{u>t\}$, that is, $\{u^*>t\}=B_{\rho(t)}$ with the same radius $\rho(t)$. (If $\{u>t\}$ is null, both sets are null and both perimeters vanish.) Applying \eqref{eq:oneball} for each such $t$ and integrating,
\[
  [u]_{0,\varepsilon}=\int_0^\infty P_b\big(B_{\rho(t)}(h(t))\big)\,dt\ \ge\ \int_0^\infty P_b\big(B_{\rho(t)}\big)\,dt=\int_0^\infty P_b(\{u^*>t\})\,dt=[u^*]_{0,\varepsilon}. \qedhere
\]
\end{proof}

\section{Proof of Theorem  \ref{Main result}}\label{proof of main result}

In this section, we prove Theorem \ref{Main result}. The proof involves a dyadic decomposition of the space based on the level sets of a measurable function, which allows us to relate the corresponding integrals to sums involving the measures of these sets. We then combine the resulting estimates with a suitable interpolation argument to obtain the desired inequality.

\smallskip

One has, for  $t_{1} \geq 0, ~t_{2} \geq 0$ with  $t_{1}+ t_{2} \geq 1$, and for  $x_{k} \geq 0$ and  $y_{k} \geq 0$,
\begin{equation}\label{ineq}
        \sum_{k \in \mathbb{Z}} x^{t_{1}}_{k} y^{t_{2}}_{k} \leq \left( \sum_{k \in \mathbb{Z}} x_{k} \right)^{t_{1}}  \left( \sum_{k \in \mathbb{Z}} y_{k} \right)^{t_{2}} \ .
    \end{equation}
    
\smallskip

The following lemma establishes a relation between the integration of a nonnegative, radially nonincreasing function. Specifically, it establishes a link between the integration over an measurable set  $\Omega$ and that over an open ball centered at the origin with the same measure.

\begin{lemma}\label{symm}
    Let  $\Omega \subset \mathbb{R}^d $ be a measurable set with $|\Omega|< \infty$, and let  $u: \mathbb{R}^{d}  \to \mathbb{R}$ be a nonnegative, radially nonincreasing function. Then,
    \begin{equation}
        \int_{\Omega} u(x)  dx \leq \int_{B(0,r)} u(x)  dx  ,
    \end{equation}
    where  $B(0,r)$ is the open ball centered at the origin and radius  $r>0$, and having the same measure as  $\Omega$.
\end{lemma}

\begin{proof}
If $\Omega=B(0,r)$ for some $r>0$, then the inequality in the statement of the lemma is an equality. Therefore, we may assume that $\Omega$ is not a ball centred at the origin with $|\Omega|>0$. Define  $\Omega' := B(0,r) \cap \Omega$. Then for all  $x \in B(0, r) \backslash \Omega'$, we have
    \begin{equation*}
        u(y) \leq u(x), \hspace{3mm} \ \forall \ y \in \Omega \backslash \Omega'.
    \end{equation*}
    Therefore, we have
    \begin{equation*}
        \int_{\Omega \backslash \Omega'} u(y) dy \leq |\Omega \backslash \Omega'| u(x).
    \end{equation*}
    Again, integrating the above inequality on  $B(0,r) \backslash \Omega'$, we obtain
    \begin{equation*}
    |B(0,r) \backslash \Omega'| \int_{\Omega \backslash \Omega'} u(y) dy \leq |\Omega \backslash \Omega'| \int_{B(0,r) \backslash \Omega'} u(x) dx.
\end{equation*}
Since  $|B(0,r)|= |\Omega|$, we have
\begin{equation*}
   |B(0,r) \backslash \Omega'| = |B(0,r)| - |\Omega'| = |\Omega| - |\Omega'| = |\Omega \backslash \Omega'|. 
\end{equation*}
Therefore, we have
\begin{equation*}
    \int_{\Omega \backslash \Omega'} u(x)  dx \leq \int_{B(0,r) \backslash \Omega'} u(x)  dx .
\end{equation*}
Hence,
\begin{equation*}
      \int_{\Omega} u(x)  dx \leq \int_{B(0,r)} u(x)  dx.
\end{equation*}
This proves the lemma.
\end{proof}
\smallskip

\begin{remark}
  \normalfont  It is obvious and well known that the above lemma  follows from a standard symmetrization argument  \cite[Section 1.3]{kesavan2006symmetrization}  of the type 
$$ \int_{\Omega} f(x)g(x) dx \leq \int_{\Omega^*} f^{*}(x)g^{*}(x) dx,$$
where $\Omega^*$ is a ball centered at origin and having same measure as  $\Omega$, and  $f^*, ~ g^*$ denotes the Schwarz symmetrization of  $f$ and  $g$ respectively. To make the presentation of this article self-contained, we present a proof of it.  
\end{remark}

We next establish a lemma that provides a decomposition of the space into disjoint sets according to the values of a measurable function, together with an estimate relating the integral of the function to the measures of these sets.
\begin{lemma}\label{ineq lemma}
    Let  $d \geq 1, ~\tau > 0$ and  $\gamma \leq 0$ such that  $\gamma \tau + d >0$. Let  $u: \mathbb{R}^d \to \mathbb{R}$ be a measurable function with compact support.  Then there exists a constant  $C= C(d, \tau, \gamma)>0$ such that
    \begin{equation}
         \int_{\mathbb{R}^d} |x|^{\gamma \tau} |u(x)|^{\tau}  dx  \leq C  \sum_{k \in \mathbb{Z}} 2^{k \tau} d^{ \frac{\gamma \tau +d}{ d}}_{k}  ,
    \end{equation}
where $d_{k} = |\{ x \ : \ 2^{k}< |u(x)| \leq 2^{k+1} \}|$.
\end{lemma}

\begin{proof}
   Let  $u$ be a measurable function on  $\mathbb{R}^d$, and for each  $k \in \mathbb{Z}$, set
    \begin{equation*}
        D_{k}  :=  \{ x \ : \ 2^{k}< |u(x)| \leq 2^{k+1} \},
    \end{equation*}
    with  $d_{k} = |D_{k}|$. Now consider
    \begin{align*}\label{eqn010}
        \int_{\mathbb{R}^d} |x|^{\gamma \tau} |u(x)|^{\tau}  dx  =    \sum_{k \in \mathbb{Z}} \int_{D_{k}}|x|^{\gamma \tau} |u(x)|^{\tau}  dx &
        \leq  \sum_{k \in \mathbb{Z}} 2^{(k+1) \tau} \int_{D_{k}} |x|^{\gamma \tau}  dx  \\ &
        = 2^{\tau}  \sum_{k \in \mathbb{Z}} 2^{k \tau} \int_{D_{k}}|x|^{\gamma \tau}  dx.
    \end{align*}
   Now, using Lemma  \ref{symm} with  $|B(0,r_{k})|= |D_{k}|$, we have
    \begin{equation}\label{eqn01}
      \int_{\mathbb{R}^d} |x|^{\gamma \tau} |u(x)|^{\tau}  dx 
      \leq 2^{\tau} \sum_{k \in \mathbb{Z}} 2^{k \tau} \int_{B(0, r_{k})}|x|^{\gamma \tau} dx.
      \end{equation}
     Then,  $r_{k} = C d^{\frac{1}{d}}_{k}$ for some constant  $C=C(d)>0$. Therefore, we have
\begin{equation}\label{eqn02}
     \sum_{k \in \mathbb{Z}} 2^{k \tau} \int_{B(0, r_{k})}|x|^{\gamma \tau}  dx  = C  \sum_{k \in \mathbb{Z}} 2^{k \tau} \int_{0}^{r_{k}} r^{\gamma \tau +d-1}  dr  = C  \sum_{k \in \mathbb{Z}} 2^{k \tau}  r^{\gamma \tau +d}_{k}  = C  \sum_{k \in \mathbb{Z}} 2^{k \tau} d^{\frac{\gamma \tau +d}{d}}_{k} .
\end{equation}
Combining  \eqref{eqn01} and  \eqref{eqn02}, we obtain   
    \begin{equation*}
        \int_{\mathbb{R}^d} |x|^{\gamma \tau} |u(x)|^{\tau}  dx  \leq C  \sum_{k \in \mathbb{Z}} 2^{k \tau} d^{\frac{\gamma \tau +d}{d}}_{k} .
    \end{equation*}
    This finishes the proof of the lemma.
\end{proof}
\smallskip

The next lemma establishes a reverse type of relation for a nonnegative, radially nonincreasing function as obtained in the previous lemma.

\begin{lemma}\label{beta dependance}
Let  $ q > 0, ~ \beta \in \mathbb{R}$ and  $\beta q+d > 0$. Let  $u: \mathbb{R}^d \to \mathbb{R}$ be a nonnegative, radially nonincreasing function with compact support. Then there exist a constant  $C=C(d,q, \beta)>0$ such that
\begin{equation}\label{eqn009}
         \int_{\mathbb{R}^d} |x|^{\beta q} |u(x)|^{q}  dx \geq C \sum_{k \in \mathbb{Z}} 2^{k q} d^{ \frac{\beta q +d}{d}}_{k} ,
    \end{equation}
where $d_{k} = |\{ x \ : \ 2^{k}< |u(x)| \leq 2^{k+1} \}|$.  If  $\beta = 0$, then  \eqref{eqn009} holds true for any measurable function.
\end{lemma}
\begin{proof}
Let $A_{k} = \{x: |u(x)|>2^{k} \}$, $D_k=A_k\setminus A_{k+1}$,
and $d_k=|D_k|$ for each $k\in\mathbb{Z}$. Since $u$ is a nonnegative, radially nonincreasing function, there exists $r_k>0$ such that $A_k=B(0,r_k)$ up to a set of measure zero. Set $a_k=|A_k|$. Then
\begin{align*}
a_k =|B(0,r_k)|
=\frac{|\mathbb{S}^{d-1}|}{d}r_k^d.
\end{align*}
Define
\begin{equation*}
W(E):=\int_E |x|^{\beta q}\,dx.
\end{equation*}
Since $\beta q+d>0$, we obtain
\begin{align*}
W(A_k)=\int_{A_k}|x|^{\beta q}\,dx =|\mathbb{S}^{d-1}|
\int_0^{r_k}r^{\beta q+d-1}\,dr =\frac{|\mathbb{S}^{d-1}|}{\beta q+d}r_k^{\beta q+d}=\frac{|\mathbb{S}^{d-1}|}{\beta q+d}
\left(\frac{d}{|\mathbb{S}^{d-1}|}\right)^{\frac{\beta q+d}{d}}
a_k^{\frac{\beta q+d}{d}}.
\end{align*}
Since $a_k\geq d_k$ and $\beta q+d>0$, it follows that
\begin{align*}
W(A_k)\geq
\frac{|\mathbb{S}^{d-1}|}{\beta q+d}
\left(\frac{d}{|\mathbb{S}^{d-1}|}\right)^{\frac{\beta q+d}{d}}
d_k^{\frac{\beta q+d}{d}}.
\end{align*}
Since $u(x)>2^k$ on $D_k$, and $u$ has compact support,  we have
\begin{align*}
\int_{\mathbb{R}^d}|x|^{\beta q}|u(x)|^q\,dx
&\geq
\sum_{\substack{k\in\mathbb Z\\D_k\neq\emptyset}}2^{kq}W(D_k) \\
&=\sum_{\substack{k\in\mathbb Z\\D_k\neq\emptyset}}2^{kq}
\bigl(W(A_k)-W(A_{k+1})\bigr) \\ & =\sum_{\substack{k\in\mathbb Z\\D_k\neq\emptyset}}2^{kq}
W(A_k)- \sum_{\substack{k\in\mathbb Z\\D_k\neq\emptyset}}2^{kq} W(A_{k+1})
\\ &=(1-2^{-q})\sum_{\substack{k\in\mathbb Z\\D_k\neq\emptyset}}2^{kq}W(A_k) \geq C(d,q,\beta)
\sum_{\substack{k\in\mathbb Z\\D_k\neq\emptyset}}2^{kq}
d_k^{\frac{\beta q+d}{d}},
\end{align*}
where
\begin{equation*}
C(d,q,\beta)
=(1-2^{-q})
\frac{|\mathbb{S}^{d-1}|}{\beta q+d}
\left(\frac{d}{|\mathbb{S}^{d-1}|}\right)^{\frac{\beta q+d}{d}}
>0.
\end{equation*}

When $\beta=0$, the same argument applies to any measurable
function $u$, since $W(E)=|E|$ and the argument does not require
the sets $A_k$ to be balls. This completes the proof.
\end{proof}

\smallskip
 
We will now establish a lemma that gives a connection between weighted Gagliardo seminorm and the measure of a disjoint union of sets obtained from the range of nonnegative, radially nonincreasing function. A similar type of lemma is available in  \cite{di2012hitchhikers} as a combination of \cite[Lemma $6.2$]{di2012hitchhikers} and \cite[Lemma $6.3$]{di2012hitchhikers}.

\begin{lemma}\label{fractional ineq}
Let $p\ge 1$ and $s\in(0,1)$. Let $\alpha_1,\alpha_2\in\R$ satisfy $\alpha_1 p>-d$,  $\alpha_2<s$, and $\frac{1}{p} + \frac{\alpha-s}{d}>0$. Let $u\in L^\infty(\R^d)$ be nonnegative, radially nonincreasing, with compact
support. For $k\in\mathbb{Z}$ set $A_k=\{u>2^k\}$, $a_k=|A_k|$, and $d_k=a_k-a_{k+1}$. Then there is $C=C(d,s,p,\alpha_1,\alpha_2)>0$ such that
\begin{equation}\label{eq:concl}
 [u]_{W^{s,p, \alpha}}^p \;\ge\; C\sum_{k\in\mathbb{Z}}2^{pk}a_k^{\frac{\alpha p }{d} + \frac{d-sp}{d}}
  \;\ge\; C\sum_{k\in\mathbb{Z}}2^{pk}d_k^{\frac{\alpha p }{d} + \frac{d-sp}{d}}.
\end{equation}
\end{lemma}

\begin{proof}
If $u=0$ a.e.\ there is nothing to prove, so assume $u\not\equiv 0$. Since $u$ is
radially nonincreasing, each $A_k$ is a centred ball (open or closed) up to a
null set; write $A_k=B(0,r_k)$ with $a_k=|B_1|r_k^d$. The sequence $(a_k)$ is
nonincreasing and bounded by $|\operatorname{supp} u|$, $a_k=0$ whenever
$2^k\ge\|u\|_\infty$, and $a_k>0$ for some $k$.

\medskip\noindent\textit{Step 1:} Put
\begin{equation*}
J_k:=\int_{A_{k+1}}\int_{A_k^c}
  \frac{|x|^{\alpha_1 p}|y|^{\alpha_2 p}}{|x-y|^{d+sp}}\,dy\,dx .  
\end{equation*}
We claim that
\begin{equation}\label{eq:layer1}
 [u]_{W^{s,p, \alpha}(\mathbb{R}^{d})}^p\;\ge\;(1-2^{-p})\sum_{k\in\mathbb{Z}}2^{pk}J_k .
\end{equation}
Fix $x,y\in\R^d$ and let
$K(x,y):=\{k\in\mathbb{Z}:\ u(x)>2^{k+1},\ u(y)\le 2^k\}$. First assume that $K(x,y) \neq \emptyset$. Then
\begin{equation}\label{eq:ind}
\sum_{k\in\mathbb{Z}}2^{pk}\,\mathbf 1_{A_{k+1}}(x)\,\mathbf 1_{A_k^c}(y)
  =\sum_{k\in K(x,y)}2^{pk},
\end{equation}
where
\begin{equation*}
    \mathbf 1_{A_{k+1}}(x)\,\mathbf 1_{A_k^c}(y)=
  \begin{cases}
    1, & k\in K(x,y),\\
    0, & k\notin K(x,y).
  \end{cases}
\end{equation*}
Multiply \eqref{eq:ind} by the nonnegative kernel
$|x|^{\alpha_1 p}|y|^{\alpha_2 p}|x-y|^{-d-sp}$, integrate over
$\R^d\times\R^d$, and exchange sum and integral by Tonelli's theorem. This gives
\begin{equation}\label{eq876}
\sum_{k\in\mathbb{Z}}2^{pk}J_k = \sum_{k \in \mathbb{Z}} 2^{pk} \int_{A_{k+1}}\int_{A_k^c}
  \frac{|x|^{\alpha_1 p}|y|^{\alpha_2 p}}{|x-y|^{d+sp}}\,dy\,dx = \int_{\mathbb{R}^{d}} \int_{\mathbb{R}^{d}} \left( \sum_{k\in K(x,y)}2^{pk} \right) \frac{|x|^{\alpha_1 p}|y|^{\alpha_2 p}}{|x-y|^{d+sp}}\,dy\,dx .
\end{equation}
If $K(x,y)=\emptyset$, the sum in \eqref{eq:ind} is $0\le |u(x)-u(y)|^p/(1-2^{-p})$. Otherwise
$K(x,y)$ is bounded above, because $2^{k+1}<u(x)\le\|u\|_\infty$; let $k_0$ be
its largest element. Then
\begin{equation*}
    u(x)-u(y)>2^{k_0+1}-2^{k_0}=2^{k_0}.
\end{equation*}
Using the above estimate, we obtain
\begin{equation*}
  \sum_{k\in K(x,y)}2^{pk}\le\sum_{k\le k_0}2^{pk}
  =\frac{2^{pk_0}}{1-2^{-p}}\le\frac{|u(x)-u(y)|^p}{1-2^{-p}} .
\end{equation*}
Multiply the above inequality by the nonnegative kernel
$|x|^{\alpha_1 p}|y|^{\alpha_2 p}|x-y|^{-d-sp}$, integrate over
$\R^d\times\R^d$, and using the definition of $[u]_{W^{s,p,\alpha}(\mathbb{R}^{d})}$, we obtain
\begin{align*}
    \int_{\mathbb{R}^{d}} \int_{\mathbb{R}^{d}} \left( \sum_{k\in K(x,y)}2^{pk} \right) \frac{|x|^{\alpha_1 p}|y|^{\alpha_2 p}}{|x-y|^{d+sp}}\,dy\,dx & \leq \frac{1}{1-2^{-p}}   \int_{\mathbb{R}^{d}} \int_{\mathbb{R}^{d}} \frac{|x|^{\alpha_{1}p}|y|^{\alpha_{2}p}|u(x)-u(y)|^{p}}{|x-y|^{d+sp}} \, dy \, dx \\ & = \frac{[u]_{W^{s,p, \alpha}(\mathbb{R}^{d})}^p}{1-2^{-p}}.
\end{align*}
Combining the above inequality with \eqref{eq876}, we get
\begin{equation*}
\sum_{k\in\mathbb{Z}}2^{pk}J_k\le\frac{[u]_{W^{s,p, \alpha}(\mathbb{R}^{d})}^p}{1-2^{-p}}, 
\end{equation*}
which is \eqref{eq:layer1}. The weight is not symmetric in $(x,y)$, but this
causes no problem. Only pairs with $x$ in the inner set and $y$ in the outer set
are counted, which can only make the left-hand side smaller.

\medskip\noindent\textit{Step 2: } If $a_{k+1}=0$ then
$J_k\ge 0$ and the term contributes nothing. Let $a_{k+1}>0$ and put
$\rho:=r_k\ge r_{k+1}>0$. Restrict the $y$-integral to the annulus
$\{\rho\le|y|<2\rho\}$, which lies in $A_k^c$ up to a null set. For
$x\in A_{k+1}$ and $y$ in the annulus,
\[
  |x-y|\le r_{k+1}+2\rho\le 3\rho,\qquad
  |y|^{\alpha_2 p}\ge\min\{1,2^{\alpha_2 p}\}\,\rho^{\alpha_2 p}.
\]
This is valid for either sign of $\alpha_2$. Hence, for every $x\in A_{k+1}$,
\begin{align*}
   \int_{A_k^c}\frac{|y|^{\alpha_2 p}}{|x-y|^{d+sp}}\,dy & \geq  \int_{\{ y : \rho \leq |y| < 2 \rho \}}\frac{|y|^{\alpha_2 p}}{|x-y|^{d+sp}}\,dy \\ &
  \ge\min\{1,2^{\alpha_2 p}\}\,\rho^{\alpha_2 p}(3\rho)^{-d-sp}(2^d-1)|B_1|\rho^d
  =c_1\rho^{(\alpha_2-s)p},
\end{align*}
where $|B_{1}| = \frac{|\mathbb{S}^{d-1}|}{d}$. Since $\alpha_1 p>-d$,
\[
  \int_{A_{k+1}}|x|^{\alpha_1 p}\,dx=\frac{|\mathbb{S}^{d-1}|}{d+\alpha_1 p}\,
  r_{k+1}^{\,d+\alpha_1 p}.
\]
Using $r_j=(a_j/|B_1|)^{1/d}$, we obtain
\begin{equation}\label{eq:Jk}
  J_k\ge c_2\,a_{k+1}^{\,1+\frac{\alpha_1 p}{d}}\,a_k^{\frac{(\alpha_2-s)p}{d}}
  \quad\text{whenever }a_{k+1}>0,
\end{equation}
with $c_2=c_2(d,s,p,\alpha_1,\alpha_2)>0$. Combining \eqref{eq:layer1} and
\eqref{eq:Jk}, and shifting $k+1\mapsto k$,
\begin{equation}\label{eq:X}
 [u]_{W^{s,p, \alpha}}^p\ge c_3X,\qquad
  X:=\sum_{k:\,a_k>0}2^{pk}a_k^{A}a_{k-1}^{B},\qquad
  A:=1+\frac{\alpha_1 p}{d}>0,\quad B:=\frac{(\alpha_2-s)p}{d}<0,
\end{equation}
where $c_3=(1-2^{-p})2^{-p}c_2$. If $a_k>0$ then $a_{k-1}\ge a_k>0$, so the
negative power $a_{k-1}^{B}$ is well defined.

\medskip\noindent\textit{Step 3: } Let
\[
  q_{*}=A+B=\frac{(\alpha_1+\alpha_2)p+d-sp}{d} = \frac{\alpha p}{d} + \frac{d-sp}{d}>0
\]
Hence
$\lambda:=q_{*}/A\in(0,1)$. For every $k$ with $a_k>0$,
\begin{equation}\label{eq:holder}
  2^{pk}a_k^{\,q_{*}}=\bigl(2^{pk}a_k^{A}a_{k-1}^{B}\bigr)^{\lambda}
  \bigl(2^{pk}a_{k-1}^{\,q_{*}}\bigr)^{1-\lambda}.
\end{equation}
Indeed, the exponent of $a_k$ is $\lambda A=q_{*}$. The exponent of $a_{k-1}$ is
\[
  \lambda B+(1-\lambda)q_{*}=\frac {q_{*}}{A}\,(B+A-q_{*})=0,
\]
and the exponent of $2^{pk}$ is $\lambda+1-\lambda=1$. Set
\[
  S:=\sum_{k\in\mathbb{Z}}2^{pk}a_k^{\,q_{*}}=\sum_{k:\,a_k>0}2^{pk}a_k^{\,q_{*}}.
\]
Then $0<S<\infty$. Positivity holds because some $a_k>0$. Finiteness holds
because $a_k\le|\operatorname{supp}u|$, $a_k=0$ for $2^k\ge\|u\|_\infty$, and
$\sum_{k\le k_1}2^{pk}<\infty$ for every $k_1$. Summing \eqref{eq:holder} over
$\{a_k>0\}$ and applying H\"older's inequality with exponents $1/\lambda$ and
$1/(1-\lambda)$,
\[
  S\le X^{\lambda}\Bigl(\sum_{k:\,a_k>0}2^{pk}a_{k-1}^{\,q_{*}}\Bigr)^{1-\lambda}
  \le X^{\lambda}\bigl(2^{p}S\bigr)^{1-\lambda},
\]
since $\sum_k 2^{pk}a_{k-1}^{\,q_{*}} \leq 2^p\sum_k 2^{pk}a_k^{\,q_{*}}$. Because
$0<S<\infty$, we may divide by $S^{1-\lambda}$:
\[
  S^{\lambda}\le 2^{p(1-\lambda)}X^{\lambda},\qquad\text{that is,}\qquad
  S\le 2^{\frac{p(1-\lambda)}{\lambda}}X.
\]
Together with \eqref{eq:X} this proves the first inequality in \eqref{eq:concl},
with $C=c_3\,2^{-p(1-\lambda)/\lambda}$. The second inequality follows from
$0\le d_k\le a_k$ and $q_{*}>0$.
\end{proof}

\begin{remark}\label{Rem: 1234}
When $\alpha_1=\alpha_2=0$, the inequality \eqref{eq:concl} becomes
\begin{equation*}
[u]_{W^{s,p}(\mathbb{R}^d)}^p \geq C\sum_{k\in\mathbb{Z}}2^{pk}a_k^{\frac{d-sp}{d}}
\geq C\sum_{k\in\mathbb{Z}}2^{pk}d_k^{\frac{d-sp}{d}},
\end{equation*}
and holds for any compactly supported
$u\in L^\infty(\mathbb{R}^d)$, without the assumption that $u$ is
radially nonincreasing. For the proof, we refer to \cite[Lemma $6.2$ and Lemma $6.3$]{di2012hitchhikers}.
\end{remark}

\begin{proof}[\textbf{Proof of Theorem  \ref{Main result}}]
Let $u \in C^{1}_{c}(\mathbb{R}^d)$ be a nonnegative, radially nonincreasing function. From Lemma \ref{ineq lemma}, we have
    \begin{equation*}
       \left( \int_{\mathbb{R}^d} |x|^{\gamma \tau} |u(x)|^{\tau} dx \right)^{\frac{1}{\tau}} \leq C \left( \sum_{k \in \mathbb{Z}} 2^{k \tau} d^{ \frac{\gamma \tau +d}{ d}}_{k} \right)^{\frac{1}{\tau}}  \  . 
    \end{equation*}
    Using the definition of  $\frac{1}{\tau} + \frac{\gamma}{d}$, we have
    \begin{align*}
            \left( \int_{\mathbb{R}^d} |x|^{\gamma \tau} |u(x)|^{\tau}  dx \right)^{\frac{1}{\tau}}  & \leq C \left( \sum_{k \in \mathbb{Z}} 2^{k \tau} d^{ a \left( \frac{1}{p} + \frac{\alpha-s}{d} \right)\tau + (1-a) \left( \frac{1}{q} + \frac{\beta}{d} \right) \tau}_{k} \right)^{\frac{1}{\tau}} \\ &
        = C \left( \sum_{k \in \mathbb{Z}} \left( 2^{kp} d^{\frac{\alpha p}{d} + \frac{d-sp}{d} }_{k} \right)^{\frac{ a \tau}{p}} \left( 2^{kq} d^{\frac{\beta q + d}{d}}_{k} \right)^{\frac{(1-a)\tau}{q}} \right)^{\frac{1}{\tau}} .
    \end{align*}
    Applying  \eqref{ineq}  with  $t_{1} =  \frac{a \tau}{p}$ and  $t_{2} =  \frac{(1-a) \tau}{q}$, we obtain
     \begin{equation}\label{eq11}
        \left( \int_{\mathbb{R}^d} |x|^{\gamma \tau} |u(x)|^{\tau}  dx \right)^{\frac{1}{\tau}} \leq C \left( \sum_{k \in \mathbb{Z}} 2^{kp} d^{\frac{\alpha p}{d} + \frac{d-sp}{d} }_{k} \right)^{\frac{a}{p}} \left( \sum_{k \in \mathbb{Z}} 2^{kq} d^{\frac{\beta q + d}{d}}_{k} \right)^{\frac{1-a}{q}}  .
    \end{equation}
    Therefore, using Lemma  \ref{beta dependance} and Lemma  \ref{fractional ineq}, and  the inequality  \eqref{eq11}, we have
    \begin{equation*}
         \| |x|^{\gamma}u \|_{L^{\tau}(\mathbb{R}^d)} \leq C [u]^{a}_{W^{s,p, \alpha}(\mathbb{R}^d)} \| |x|^{\beta} u \|^{(1-a)}_{L^{q}(\mathbb{R}^d)} \ .
    \end{equation*}
    This proves the first part of Theorem  \ref{Main result}.
    
\smallskip

Now, assume  $\alpha_{1}=\alpha_{2}=0$ and  $\beta=0$ and let  $u \in C^{1}_{c}(\mathbb{R}^d)$. From Lemma  \ref{ineq lemma} , we have
    \begin{equation*}
       \left( \int_{\mathbb{R}^d} |x|^{\gamma \tau} |u(x)|^{\tau} dx \right)^{\frac{1}{\tau}} \leq C \left( \sum_{k \in \mathbb{Z}} 2^{k \tau} d^{ \frac{\gamma \tau +d}{ d}}_{k} \right)^{\frac{1}{\tau}}  \  . 
    \end{equation*}
    Using the definition of  $\frac{1}{\tau} + \frac{\gamma}{d}$, we have 
    \begin{align*}
            \left( \int_{\mathbb{R}^d} |x|^{\gamma \tau} |u(x)|^{\tau}  dx \right)^{\frac{1}{\tau}} & \leq C \left( \sum_{k \in \mathbb{Z}} 2^{k \tau} d^{ a \left( \frac{1}{p} - \frac{s}{d} \right)\tau +   \frac{(1-a)\tau}{q}}_{k} \right)^{\frac{1}{\tau}} \\ &
        = C \left( \sum_{k \in \mathbb{Z}} \left( 2^{kp} d^{\frac{d-sp}{d}  }_{k} \right)^{\frac{ a \tau}{p}} \left( 2^{kq} d_{k} \right)^{\frac{(1-a)\tau}{q}} \right)^{\frac{1}{\tau}} .
    \end{align*}
Applying  \eqref{ineq}  with  $t_{1} =  \frac{a \tau}{p}$ and  $t_{2} =  \frac{(1-a) \tau}{q}$, we obtain
    \begin{equation}\label{eq4}
        \left( \int_{\mathbb{R}^d} |x|^{\gamma \tau} |u(x)|^{\tau}  dx \right)^{\frac{1}{\tau}} \leq C \left( \sum_{k \in \mathbb{Z}} 2^{kp} d^{\frac{d-sp}{d}}_{k} \right)^{\frac{a}{p}} \left( \sum_{k \in \mathbb{Z}} 2^{kq} d_{k} \right)^{\frac{(1-a)}{q}}  .
    \end{equation}
    In the above step we have used the assumption $\tau \geq \max\{p,q\}$ in order to ensure $t_1+t_2 \geq 1$.
Therefore, using Lemma  \ref{beta dependance} with  $\beta=0$ and Remark \ref{Rem: 1234} with  $\alpha_{1}=\alpha_{2}=0$, and inequality \eqref{eq4}, we have
    \begin{equation*}
        \| |x|^{\gamma}u \|_{L^{\tau}(\mathbb{R}^d)} \leq C [u]^{a}_{W^{s,p}(\mathbb{R}^d)} \| u \|^{(1-a)}_{L^{q}(\mathbb{R}^d)}  .
    \end{equation*}
    \end{proof}

\begin{remark}
   Note that in order to get   \eqref{eq4} it is sufficient to assume the weaker condition  $\frac{a\tau}{p}+\frac{(1-a)\tau}{q}\geq 1$ than assuming  $\tau \geq \max\{p,q\}$.
\end{remark}

\bigskip

\textbf{Acknowledgement:} P. Roy is supported by the ARG MATRICS Grant  \\(ANRF/ARGM/2025/002269/MTR) of the Anusandhan National Research Foundation (ANRF). V. Sahu is supported by the National Postdoctoral Fellowship (PDF/2025/004611) of the Anusandhan National Research Foundation (ANRF). P. Roy would like to thank Adimurthi for helpful discussions on the subject and for pointing out the main idea.

\smallskip

\textbf{Use of AI tools: }  \textit{Claude opus 5.5 medium } was used for this work. However, the authors take full responsibility for the mathematical reasoning.


\begin{thebibliography}{10}

\bibitem{Abdellaoui2017}
B.~Abdellaoui and R.~Bentifour, \emph{Caffarelli-{K}ohn-{N}irenberg type inequalities of fractional order with applications}, J. Funct. Anal. \textbf{272} (2017), no.~10, 3998--4029.

\bibitem{abdellaoui2005some}
B.~Abdellaoui, E.~Colorado, and I.~Peral, \emph{Some improved {C}affarelli-{K}ohn-{N}irenberg inequalities}, Calc. Var. Partial Differential Equations \textbf{23} (2005), no.~3, 327--345.

\bibitem{adi2002}
Adimurthi, N.~Chaudhuri, and M.~Ramaswamy, \emph{An improved {H}ardy-{S}obolev inequality and its application}, Proc. Amer. Math. Soc. \textbf{130} (2002), no.~2, 489--505.

\bibitem{ADICCM2026}
Adimurthi, P.~Jana, and P.~Roy, \emph{Boundary fractional {H}ardy's inequality in dimension one: the critical case}, Commun. Contemp. Math. \textbf{28} (2026), no.~4, Paper No. 2550051, 14.

\bibitem{adiarka2018}
Adimurthi and A.~Mallick, \emph{A {H}ardy type inequality on fractional order {S}obolev spaces on the {H}eisenberg group}, Ann. Sc. Norm. Super. Pisa Cl. Sci. (5) \textbf{18} (2018), no.~3, 917--949.

\bibitem{AdiJFA2026}
Adimurthi, P.~Roy, and V.~Sahu, \emph{Fractional boundary {H}ardy inequality for the critical cases}, J. Funct. Anal. \textbf{290} (2026), no.~8, Paper No. 111351, 63.

\bibitem{AdiCVPDE2026}
Adimurthi, P.~Roy, and V.~Sahu, \emph{Fractional {H}ardy inequality with singularity on submanifold}, Calc. Var. Partial Differential Equations \textbf{65} (2026), no.~2, Paper No. 62, 40.

\bibitem{Adi2025}
Adimurthi, P.~Roy, and V.~Sahu, \emph{Boundary {H}ardy inequality on functions of bounded variation}, Ann. Sc. Norm. Super. Pisa Cl. Sci. (5) (Mar. 2025), p. 26. \href{https://journals.sns.it/index.php/annaliscienze/article/view/6798}{https://journals.sns.it/index.php/annaliscienze/article/view/6798}.

\bibitem{Almgren1989}
F.~J. Almgren, Jr. and E.~H. Lieb, \emph{Symmetric decreasing rearrangement is sometimes continuous}, J. Amer. Math. Soc. \textbf{2} (1989), no.~4, 683--773.

\bibitem{Ao2022}
W.~Ao, A.~DelaTorre, and M.~del~M. Gonz\'alez, \emph{Symmetry and symmetry breaking for the fractional {C}affarelli-{K}ohn-{N}irenberg inequality}, J. Funct. Anal. \textbf{282} (2022), no.~11, Paper No. 109438, 58.

\bibitem{BanerjeeGangulySahu}
A.~Banerjee, D.~Ganguly, and V.~Sahu, \emph{Quantitative stability for fractional {H}ardy inequalities: Rearrangement-free techniques and {E}mden-{F}owler analysis}, Preprint, arXiv:2605.15748 (2026), 44 pages.

\bibitem{bartsch2007existence}
T.~Bartsch, S.~Peng, and Z.~Zhang, \emph{Existence and non-existence of solutions to elliptic equations related to the Caffarelli-Kohn-Nirenberg inequalities}, Calculus of Variations and Partial Differential Equations \textbf{30} (2007), no.~1, 113--136.

\bibitem{bogdan2011dyda}
K.~Bogdan and B.~Dyda, \emph{The best constant in a fractional Hardy inequality}, Math. Nachr. \textbf{284} (2011), no.~5-6, 629--638.

\bibitem{Brezis1997}
H.~Brezis and M.~Marcus, \emph{Hardy's inequalities revisited}, Ann. Scuola Norm. Sup. Pisa Cl. Sci. (4) \textbf{25} (1997), no.~1-2, 217--237, Dedicated to Ennio De Giorgi.

\bibitem{Brezis2000}
H.~Brezis, M.~Marcus, and I.~Shafrir, \emph{Extremal functions for {H}ardy's inequality with weight}, J. Funct. Anal. \textbf{171} (2000), no.~1, 177--191.

\bibitem{Caffarelli1984}
L.~Caffarelli, R.~Kohn, and L.~Nirenberg, \emph{First order interpolation inequalities with weights}, Compositio Math. \textbf{53} (1984), no.~3, 259--275.

\bibitem{catrina2000caffarelli}
F.~Catrina and Z.-Q. Wang, \emph{On the {C}affarelli-{K}ohn-{N}irenberg inequalities}, C. R. Acad. Sci. Paris S\'er. I Math. \textbf{330} (2000), no.~6, 437--442.

\bibitem{catrina2001caffarelli}
F.~Catrina and Z.-Q. Wang, \emph{On the {C}affarelli-{K}ohn-{N}irenberg inequalities: sharp constants, existence (and nonexistence), and symmetry of extremal functions}, Comm. Pure Appl. Math. \textbf{54} (2001), no.~2, 229--258.

\bibitem{Chen2021}
L.~Chen, G.~Lu, and C.~Zhang, \emph{Maximizers for fractional {C}affarelli-{K}ohn-{N}irenberg and {T}rudinger-{M}oser inequalities on the fractional {S}obolev spaces}, J. Geom. Anal. \textbf{31} (2021), no.~4, 3556--3582.

\bibitem{chern2010minimizers}
J.-L. Chern and C.-S. Lin, \emph{Minimizers of {C}affarelli-{K}ohn-{N}irenberg inequalities with the singularity on the boundary}, Arch. Ration. Mech. Anal. \textbf{197} (2010), no.~2, 401--432. \MR{2660516}

\bibitem{Gyula2026}
G.~Csat{\'o} and P.~Roy, \emph{A fractional {H}ardy-{S}obolev inequality of {M}ichael--{S}imon type on convex hypersurfaces}, Math. Ann. \textbf{395} (2026), no.~1, 14.

\bibitem{DeNitti2025}
N.~De~Nitti, F.~Glaudo, and T.~K\"onig, \emph{Non-degeneracy, stability, and symmetry for the fractional {C}affarelli-{K}ohn-{N}irenberg inequality}, Comm. Partial Differential Equations \textbf{50} (2025), no.~8, 984--1028.

\bibitem{di2012hitchhikers}
E.~Di~Nezza, G.~Palatucci, and E.~Valdinoci, \emph{Hitchhiker's guide to the fractional {S}obolev spaces}, Bull. Sci. Math. \textbf{136} (2012), no.~5, 521--573.

\bibitem{Serena2015}
S.~Dipierro and E.~Valdinoci, \emph{A density property for fractional weighted {S}obolev spaces}, Atti Accad. Naz. Lincei Rend. Lincei Mat. Appl. \textbf{26} (2015), no.~4, 397--422.

\bibitem{dolbeaulty2009symmetry}
J.~Dolbeault, M.~J. Esteban, M.~Loss, and G.~Tarantello, \emph{On the symmetry of extremals for the {C}affarelli-{K}ohn-{N}irenberg inequalities}, Adv. Nonlinear Stud. \textbf{9} (2009), no.~4, 713--726.

\bibitem{dyda2004fractional}
B.~Dyda, \emph{A fractional order {H}ardy inequality}, Illinois J. Math. \textbf{48} (2004), no.~2, 575--588.

\bibitem{Dyda2024JFA}
B.~Dyda and M.~Kijaczko, \emph{Sharp fractional {H}ardy inequalities with a remainder for {$1 < p < 2$}}, J. Funct. Anal. \textbf{286} (2024), no.~9, Paper No. 110373, 19.

\bibitem{dyda2024}
B.~Dyda and M.~Kijaczko, \emph{Sharp weighted fractional {H}ardy inequalities}, Studia Math. \textbf{274} (2024), no.~2, 153--171.

\bibitem{dyda2015}
B.~Dyda and A.~V. V\"{a}h\"{a}kangas, \emph{Characterizations for fractional Hardy inequality}, Adv. Calc. Var. \textbf{8} (2015), no.~2, 173--182.

\bibitem{frank2008non}
R.~L. Frank and R.~Seiringer, \emph{Non-linear ground state representations and sharp {H}ardy inequalities}, J. Funct. Anal. \textbf{255} (2008), no.~12, 3407--3430.

\bibitem{Gonzalez2025}
M.~del~M. Gonz\'alez, A.~Hyder, and M.~S\'aez, \emph{The limiting case of the fractional {C}affarelli-{K}ohn-{N}irenberg inequality in dimension one}, Int. Math. Res. Not. IMRN (2025), no.~18, Paper No. rnaf270, 31.

\bibitem{kesavan2006symmetrization}
S.~Kesavan, \emph{Symmetrization \& applications}, Series in Analysis, vol.~3, World Scientific Publishing Co. Pte. Ltd., Hackensack, NJ, 2006.

\bibitem{KijaczkoSNS}
M.~Kijaczko, \emph{Fractional {S}obolev spaces with power weights}, Ann. Sc. Norm. Super. Pisa Cl. Sci. (5) \textbf{24} (2023), no.~3, 1741--1764.

\bibitem{VivekQS126}
M.~Kijaczko and V.~Sahu, \emph{Logarithmic fractional {H}ardy inequalities}, Preprint, arXiv:2610.38420 (2026), 22 pages.

\bibitem{VivekCCM2026}
M.~Kijaczko and V.~Sahu, \emph{Weighted fractional {H}ardy-{S}obolev and {H}ardy-{S}obolev-{M}az'ya inequalities with singularities on flat submanifold}, Commun. Contemp. Math. (2026), Paper No. 2650016, 28.

\bibitem{kristaly2013caffarelli}
A.~Krist\'aly and S.~Ohta, \emph{Caffarelli-{K}ohn-{N}irenberg inequality on metric measure spaces with applications}, Math. Ann. \textbf{357} (2013), no.~2, 711--726. \MR{3096522}

\bibitem{loss2010hardy}
M.~Loss and C.~Sloane, \emph{Hardy inequalities for fractional integrals on general domains}, J. Funct. Anal. \textbf{259} (2010), no.~6, 1369--1379.

\bibitem{mallick2022caffarelli}
A.~Mallick and H.-M. Nguyen, \emph{The {C}affarelli-{K}ohn-{N}irenberg inequalities for radial functions}, C. R. Math. Acad. Sci. Paris \textbf{361} (2023), 1175--1189. \MR{4659077}

\bibitem{Mazya2002}
V.~Maz'ya and T.~Shaposhnikova, \emph{On the {B}ourgain, {B}rezis, and {M}ironescu theorem concerning limiting embeddings of fractional {S}obolev spaces}, J. Funct. Anal. \textbf{195} (2002), no.~2, 230--238.

\bibitem{radulescu2011}
M.~Mihăilescu, V.~Rădulescu, and D.~Stancu-Dumitru, \emph{A {C}affarelli-{K}ohn-{N}irenberg-type inequality with variable exponent and applications to {PDE}s}, Complex Var. Elliptic Equ. \textbf{56} (2011), no.~7-9, 659--669.

\bibitem{Nguyen2018}
H.-M. Nguyen and M.~Squassina, \emph{Fractional {C}affarelli-{K}ohn-{N}irenberg inequalities}, J. Funct. Anal. \textbf{274} (2018), no.~9, 2661--2672.

\bibitem{nguyen2019hardy}
H.-M. Nguyen and M.~Squassina, \emph{On {H}ardy and {C}affarelli-{K}ohn-{N}irenberg inequalities}, J. Anal. Math. \textbf{139} (2019), no.~2, 773--797.

\bibitem{rawat2026}
R.~Rawat, H.~Roy, and P.~Roy, \emph{Fractional {C}affarelli-{K}ohn-{N}irenberg type inequalities on the {H}eisenberg group}, J. Differential Equations \textbf{457} (2026), Paper No. 114030, 34.

\bibitem{VivekJMAA2025}
V.~Sahu, \emph{Weighted fractional {H}ardy inequalities with singularity on any flat submanifold}, J. Math. Anal. Appl. \textbf{546} (2025), no.~2, Paper No. 129227, 16.

\bibitem{VivekQS26}
V.~Sahu, \emph{Quantitative stability for weighted {H}ardy inequalities: {L}ocal and nonlocal cases}, Preprint, arXiv:2609.25952 (2026), 46 pages.

\bibitem{wang2003caffarelli}
Z.-Q. Wang and M.~Willem, \emph{Caffarelli-{K}ohn-{N}irenberg inequalities with remainder terms}, J. Funct. Anal. \textbf{203} (2003), no.~2, 550--568.

\bibitem{wei2022stability}
J.~Wei and Y.~Wu, \emph{On the stability of the {C}affarelli-{K}ohn-{N}irenberg inequality}, Math. Ann. \textbf{384} (2022), no.~3-4, 1509--1546.

\end{thebibliography}

\end{document}